\documentclass{ieeeaccess}
\usepackage{amssymb}
\usepackage{pifont}
\newcommand{\xmark}{\ding{55}}%
\usepackage{hyperref}
\usepackage{cite}
\usepackage{subcaption}
\usepackage{amsmath,amssymb,amsfonts}
\usepackage{subcaption}
\usepackage[labelformat=parens,labelsep=quad,skip=3pt]{caption}
\usepackage{algorithmic}
\usepackage{graphicx}
\usepackage{textcomp}

\newtheorem{thm}{Theorem}[section]

\newtheorem{rem}[thm]{\bf{Remark}}
\usepackage{amsmath, amscd,amsfonts, amssymb }
\usepackage{times}
\usepackage{extarrows}
\usepackage{stmaryrd}
\usepackage{textcomp}
\usepackage[mathscr]{euscript}
\usepackage{algorithm2e}
\usepackage{colortbl}
\usepackage{color,soul}
\usepackage{newunicodechar}
\newcommand{\cA}{\mathbf{A}}

\newcommand{\cQ}{{\mathbf Q}}
\newcommand{\cB}{{\mathbf B}}
\newcommand{\cX}{{\mathbf X}}

\newcommand{\cx}{{\mathbf x}}
\def\BibTeX{{\rm B\kern-.05em{\sc i\kern-.025em b}\kern-.08em
    T\kern-.1667em\lower.7ex\hbox{E}\kern-.125emX}}
\begin{document}
\history{}
\doi{}

\title{Randomized Algorithms for Computation of Tucker Decomposition and Higher Order SVD (HOSVD)}
\author{
\uppercase{Salman AHmadi-Asl}$^*$\authorrefmark{1},
\uppercase{Stanislav Abukhovich}$^*$\authorrefmark{1},
\uppercase{Maame G. Asante-Mensah}$^*$\authorrefmark{1},
\uppercase{Andrzej Cichocki\authorrefmark{1}}, \IEEEmembership{Fellow, IEEE},
\uppercase{Anh Huy Phan}\authorrefmark{1},~\IEEEmembership{Member,~IEEE}
\uppercase{Tohishisa Tanaka}\authorrefmark{2},~\IEEEmembership{Senior Member,~IEEE}
and
\uppercase{Ivan Oseledets}\authorrefmark{1}
}
\address[1]{Skolkovo Institute of Science and Technology (SKOLTECH), CDISE, Moscow, Russia (e-mail: s.asl@skoltech.ru)}
\address[2]{Tokyo University of Agriculture and Technology, Tokyo, Japan\\
$^*$Contributed equally.}
\tfootnote{This work was partially supported by Mega Grant project (14.756.31.0001).}

\markboth
{Author \headeretal: Preparation of Papers for IEEE TRANSACTIONS and JOURNALS}
{Author \headeretal: Preparation of Papers for IEEE TRANSACTIONS and JOURNALS}


\begin{abstract}
Big data analysis has become a crucial part of new emerging technologies such as the internet of things, cyber-physical analysis, deep learning, anomaly detection, etc. Among many other techniques, dimensionality reduction plays a key role in such analyses and facilitates feature selection and feature extraction. Randomized algorithms are efficient tools for handling big data tensors. They accelerate decomposing large-scale data tensors by reducing the computational complexity of deterministic algorithms and the communication among different levels of the memory hierarchy, which is the main bottleneck in modern computing environments and architectures. In this paper, we review recent advances in randomization for the computation of Tucker decomposition and Higher Order SVD (HOSVD). We discuss random projection and sampling approaches, single-pass, and multi-pass randomized algorithms, and how to utilize them in the computation of the Tucker decomposition and the HOSVD. Simulations on synthetic and real datasets are provided to compare the performance of some of the best and most promising algorithms.
\end{abstract}


\titlepgskip=-15pt

\maketitle

\section{Introduction}\label{sec:introduction}
\begin{figure}
\begin{center}
\includegraphics[width=8 cm,height=4 cm]{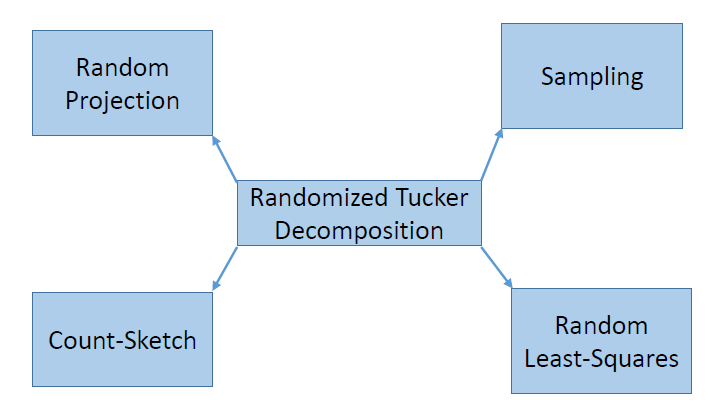}\\
\caption{{\small Taxonomy of randomized algorithms for low-rank Tucker approximation.}}\label{Category2}
\end{center}
\end{figure}

Real world data often are multidimensional and naturally are represented as tensors or multidimensional (multiway) arrays. It is crucial to preserve the multidimensional structure of the data tensors in order to extract meaningful latent variables and reveal the hidden structures of the data tensors. Tucker decomposition \cite{tucker1964extension,tucker1966some,tucker1963implications} is a natural generalization of the SVD to higher-order data tensors and has found various applications such as reducing the number of parameters in deep neural networks \cite{kim2015compression, gusak2019one}, handwritten digit classification \cite{savas2007handwritten}, computer vision \cite{vasilescu2002multilinear}, recommender systems \cite{rafailidis2012tfc, nanopoulos2009musicbox, sun2005cubesvd, symeonidis2015clusthosvd}, signal processing \cite{de2004dimensionality, cichocki2015tensor, sidiropoulos2017tensor}, etc. Deterministic algorithms for decomposing large-scale data tensors into the Tucker format are prohibitive and require high memory and computational complexity or only applicable for structured data tensors \cite{kolda2008scalable,badeau2008fast}. It has been proved that randomized algorithms can tackle this difficulty by exploiting only a part of the data tensors with applications in tensor regression \cite{zhang2020islet}, tensor completion \cite{feng2020orthogonal, ahmadi2020randomized} and deep learning \cite{kolbeinsson2021tensor}. Because of this property, they scale quite well to the tensor dimensions and orders. Two main features of the randomized algorithms which make them suitable approaches for handling large-scale data tensors are: 1) they are parallelizable and 2) they reduce the communication among different levels of memories, which is the main bottleneck in modern computing
environments and architectures. In particular, the second property is important for the data tensors stored out-of-cores, where, the communication cost may exceed our main computations. Due to the two mentioned benefits, in recent years, there is a growing interest in developing randomized algorithms for the Tucker decomposition and Higher Order SVD\footnote{It is also called MultiLinear SVD (MLSVD).} (HOSVD) \cite{de2000multilinear}. 
Randomized algorithms for the Tucker decomposition are mainly categorized into four groups as follows (see Figure \ref{Category2} and Table \ref{Categorytable})
\begin{itemize}
\item Random projection

\item Sampling 

\item Count-sketch

\item Random least-squares 
\end{itemize}

\begin{table}
\begin{center} 
\caption{Randomized Tucker references according to four categories.} 
 \begin{tabular}{||c c||}
 \hline
 Algorithms  & Ref  \\ [0.5ex] 
 \hline\hline
 Random Projection & \cite{che2018randomized, minster2019randomized, zhou2014decomposition, sun2018tensor, wolf2019low, batselier2018computing}  \\ 
  \hline\hline
 Sampling & \cite{drineas2007randomized, caiafa2010generalizing, sun2009multivis, oseledets2008tucker, saibaba2016hoid, tsourakakis2010mach, song2019relative,  perros2015sparse, cheng2016spals}  \\
 \hline
 Count-Sketch & \cite{malik2018low, shi2019multi, shi2019efficient} \\
 \hline
 Random least-squares  & \cite{malik2018low}  \\
 \hline
 \hline
\end{tabular}
\label{Categorytable}
\end{center}
\end{table}

Following the above categories, in this paper, we systematically review a variety of randomized algorithms for decomposing large-scale tensors in the Tucker and the HOSVD formats. 

The remainder of this paper is structured as follows. In Section \ref{Sec2}, we introduce preliminary concepts and notations used throughout the paper. In Section \ref{Sec:RandomMatrix}, the randomized algorithms for low-rank matrix approximation are presented. In Section \ref{Sec:HOSVD}, the Tucker decomposition and the HOSVD are briefly introduced. The randomized variants of the deterministic Tucker and HOSVD algorithms are presented in Section \ref{Sec:RHOSVD}. An application of the randomized Tucker algorithms is presented for fast computation of Canonical Polyadic Decomposition (CPD) in Section \ref{Sec:Appl}. Further challenges and open problems are discussed in Section \ref{Sec:FuthChall}. Extensive simulations on synthetic and real data tensors are provided in Section \ref{Sec:Sim} to compare the performance of some of the presented randomized algorithms. Conclusions are given in Section \ref{Sec:Concl}.

\section{Preliminary concepts and notations}\label{Sec2}
In this section, we present some concepts and notations used throughout the paper. 

Tensors\footnote{Tensors are also called multidimensional arrays or multi-ways data. From this perspective, vectors and matrices are first order and second order tensors, respectively.}, matrices and vectors are respectively denoted by underlined bold upper case letters, e.g., $\underline{\cX}$, bold upper case letters, e.g., ${\cX}$ and bold lower case letters, e.g., ${\cx}$. The number of indices or modes of a tensor is called order. Fibers are first order
tensors produced by fixing all modes except one, while slices are
second order tensors generated by fixing all modes except two
of them.

The symbols $\otimes,\,\,\dag,\,\,T$ denote the Kronecker product, the Moore-Penrose pseudoinverse and the transpose of matrices, respectively.
The Frobenius norm of tensors and spectral norm of matrices are respectively denoted by ${\left\| . \right\|_F}$ and ${\left\| . \right\|_2}$. In the case of vectors, ${\left\| . \right\|_2}$ denotes the Euclidean norm.

Tensors can be reshaped to matrices and this is called unfolding, matricization or flattening, and the corresponding generated matrices are called unfolding or flattening matrices. In the context of tensor computations, unfolding matrices are of interest because of their applications in computing different types of tensor decompositions where low-rank approximations of these matrices are required in each step of the underlying algorithms.
For example, it is known that the factor matrices of the HOSVD are computed via the left singular vectors of the unfolding matrices \cite{de2000multilinear}. The HOSVD is a direct generalization of the classical SVD of matrices to tensors where instead of rank, multilinear rank or Tucker rank is defined which is an $N$-tuple containing the rank of unfolding matrices along different modes. 
There are several types of the unfolding matrices including {\em $n$-unfolding} and {\em mode-$n$ unfolding} \cite{cichocki2016tensor}. The first one is used in computation of the Tucker decomposition and the HOSVD.

Given an $N$th-order tensor $\underline{\cX}\in\mathbb{R}^{{I_1}\times I_2\times\cdots\times {I_N}}$, then the $n$-unfolding matrix of the tensor $\underline{\cX}$ is denoted by ${\cX}_{(n)}\in\mathbb{R}^{{I_n}\times {I_1\cdots I_{n-1}I_{n+1}\cdots I_N}}$, whose components are defined as follows
\[{{\cX}_{\left( n \right)}}\left( {{i_n},j} \right) = \underline{\cX}\left( {{i_1},{i_2}, \ldots {i_N}} \right),
\]
where 
\[
j= 1+\sum\limits_{k = 1,k\ne n}^N {\left( {{i_k} - 1} \right){J_k}} ,\,\,\,\,\,\,{J_k} = \prod\limits_{m = 1,\,m\ne n}^{k - 1} {{I_m}}.
\]

Although the unfolding operator can be defined in a general setting (see \cite{kolda2006multilinear}), however throughout the paper, we use the $n$-unfolding.

Tensors and matrices can be multiplied in different modes which have the same dimensions (sizes). This is called tensor-matrix multiplication along mode $n$ or $n$-mode (matrix) product of a tensor by a matrix and is a generalization of the matrix-matrix multiplication. To be precise, let $\underline{\cX}\in\mathbb{R}^{I_1\times I_2\times \cdots\times I_N}$ and ${\cB}\in\mathbb{R}^{J\times I_n}$, then the tensor $\underline{\cX}$ and the matrix ${\cB}$ can be multiplied along mode $n$ as $\underline{\cX}\times_{n}{\cB}\in \mathbb{R}^{I_1\times \cdots \times I_{n-1}\times J\times I_{n+1}\times \cdots \times I_N}$ which is defined as follows 
\[
{\left( {{\underline{\cX}}{ \times _n}{\cB}} \right)_{{i_1},\ldots {i_{n - 1}},j,{i_{n + 1}}, \ldots, {i_N}}} =  \sum\limits_{{i_n=1}}^{{I_N}} {{x_{{i_1},{i_2}, \ldots ,{i_N}}}{b_{j,{i_n}}}},
\]
for $j=1,2,\ldots,J$. Assume that tensor $\underline{\cX}$ and two matrices ${\cA},{\cB}$ are of conforming sizes, then we have
\begin{equation}\label{TMM-M}
\underline{\cX}{ \times _n}{\cB}{ \times _n}{\cA}=\underline{\cX}{ \times _n}{\cA}{\cB}.
\end{equation}

If a tensor is multiplied with several matrices along all its modes except one, then its $n$-unfolding can be computed as follows
\begin{equation}\label{TuckerUnfold}
\begin{array}{l}
{\underline {\cX}} = {\underline {\mathbf S}}{ \times _1}{{\cQ}^{(1)}}{ \times _2}{{\cQ}^{(2)}} \ldots { \times _N}{{\cQ}^{(N)}}
\quad\Leftrightarrow 
\\
\\
\begin{array}{l}
{{\mathbf X}_{\left( n \right)}} = {{\mathbf Q}^{\left( n \right)}}{{\mathbf S}_{\left( n \right)}}\left( {{{\mathbf Q}^{\left( N \right)}} \otimes\cdots\otimes{{\mathbf Q}^{\left( {n + 1} \right)\,\,}}{\otimes}{\,\,{\mathbf Q}^{\left( {n - 1} \right)\,\,}}} {\otimes}\right.\\
\hspace{5.9cm}\left. { \cdots {\otimes }{\,\,{\mathbf Q}^{\left( 1 \right)\,\,}}} \right)^T.
\end{array}
\end{array}
\end{equation}

\section{Randomized low-rank matrix Algorithms}\label{Sec:RandomMatrix}
Randomized algorithms are utilized for finding low-rank approximation of large-scale data matrices in reasonable amounts of time and memory. They have found many applications in the area of machine learning and data mining, where we encounter large-scale matrices with more than hundreds of millions nonzero entries. 
Clearly, deterministic approaches are time-consuming in processing such data matrices and it is required to reduce the computational cost and also memory requirement of algorithms in order to process the data matrices more efficiently. Randomized algorithms can be used in such scenarios by finding a low approximation of unfolding matrices ${\cX}={\cX}_{(n)},\,n=1,2,\ldots,N$ using the following formula
\begin{eqnarray}\label{QBApp}
\nonumber
{\cX} &\cong & \left({\cQ}^{(1)}{\cQ}^{(1)^{\,T}}\right){\cX}\left({\cQ}^{(2)}{\cQ}^{(2)\,T}\right),\\
 &=& {{\cX}}{ \times _1}{{\cQ}^{(1)}}{\cQ}^{(1)\,T}{ \times _2}{{\cQ}^{(2)}}{\cQ}^{(2)\,T},
\end{eqnarray}
where ${{\cQ}^{(1)}}\in\mathbb{R}^{I_n\times \widehat{R}_n},\,{{\cQ}^{(2)}}\in{\mathbb{R}^{\left( {\prod\limits_{i \ne n} {{I_i}} } \right) \times \widetilde{R}_n}}$ ($\widehat{R}_n,\,\widetilde{R}_n< {\rm{rank}}\left( {\cX} \right)$) are approximate orthonormal bases for the column and row spaces of the matrix ${\cX}$. 

For an unfolding matrix ${\cX}={\cX}_{(n)}$, the dimension of the second mode is usually much larger than the first one, i.e., $\prod\limits_{i \ne n} {{I_i}} \gg I_n$, and because of this, it is often enough to make reduction in the second mode and just consider ${\cQ}={{\cQ}^{(1)}}$. 

The main steps of the randomized algorithms for low-rank matrix approximation are 
\begin{itemize}
\item {\bf Reduction.} Replacing an extremely large-scale matrix with a new one of smaller size compared to the original one capturing either the column or the row space of the original data matrix 
as much as possible. This step can be considered as a preprocessing step which makes reduction in the data.

\item {\bf SVD computation.} Applying classical deterministic algorithms, e.g., truncated SVD to a reduced data matrix and finding its low-rank approximation\footnote{In \cite{li2014large}, randomization is also used in the second step. This is considered as a two-step randomized algorithm for Nyström kernel matrix approximation. At the first step, authors use a sub-sampling approach after which they apply randomized SVD instead of deterministic SVD.}.

\item {\bf Recovering}. Recovering the SVD of the original data matrix from the SVD of the compressed one.
\end{itemize}
The first step can be done by following techniques
\begin{itemize}
\item Random Projection

\item Sampling 

\item Count-Sketch.
\end{itemize}

\subsection{Random Projection} 
The main idea of this technique is approximating the column (row) space of a given matrix by few new columns (rows) each of which is a linear combination of the columns (rows) of the original data matrix, where the coefficients of the linear combinations come from some probability distributions such as Gaussian, Bernoulli or uniform distribution\footnote{It is mentioned in \cite{halko2011finding} that the difference among different random matrices is negligible. We have also confirmed this in our simulations.}.  

This procedure is performed by multiplying a given matrix by a random matrix from the right-hand or left-hand side. It has been shown that this preserves the Euclidean distances among columns or rows approximately\footnote{Depending on which side the random matrix is multiplied.} \cite{ahfock2017statistical, johnson1984extensions}.
Let ${\cX}\in\mathbb{R}^{I\times J}$ be a given data matrix, and $R$ be a target rank. In the first stage of the random projection approach, we generate a random matrix ${\mathbf \Omega}= \left[ {{{\boldsymbol \omega}_1},{{\boldsymbol \omega}_2}, \ldots ,{{\boldsymbol \omega} _R}} \right] \in {\mathbb{R}^{J \times R}}$
whose components are independent and identically distributed (i.i.d) and taken from some probability distributions, e.g., Gaussian distribution, and then multiply the matrix ${\mathbf X}$ with the columns ${{\boldsymbol \omega}_r}$ as follows
\begin{equation}\label{proj}
{{\mathbf y}_r} = {\mathbf X}{{\boldsymbol\omega}_r}\in {\mathbb{R}^I},\quad r=1,2,\ldots,R.
\end{equation}
The compact form of relation \eqref{proj} is ${\mathbf Y} = {\mathbf X}\,{\boldsymbol\Omega},$
where ${\mathbf Y} = \left[ {{{\mathbf y}_1}, {{\mathbf y}_2}\ldots ,{{\mathbf y}_R}} \right]\in \mathbb{R}^{I \times R}$. 
The matrix $\mathbf Y$ has a smaller size than the matrix ${\mathbf X}$ and is an approximation for the range of ${\mathbf X}$. Note that the columns of the matrix $\mathbf Y$ are generally not orthogonal and in order to compute an orthogonal projection onto the range of $\mathbf Y$, it is necessary to generate a new orthonormal basis of ${\mathbf Y}$. This can be done by computing the economic QR decomposition of the matrix $\mathbf Y$ as ${\mathbf Y}={\mathbf Q}{\mathbf R}$ and using the matrix $\mathbf Q\in\mathbb{R}^{I\times R}$ in subsequent computations. 
Here, we have 
\begin{equation}\label{RQR}
{\mathbf{X}} \cong {\mathbf{Q}}{{\mathbf{Q}}^T}{\mathbf{X}},
\end{equation}
and having computed the SVD of the matrix ${\mathbf B}={\mathbf Q}^T{\mathbf X}\in\mathbb{R}^{R\times J}$, as ${\mathbf B} = {\mathbf U}{\mathbf S} {{\mathbf V}^T},$ the SVD of the matrix ${\mathbf X}$ is recovered as follows
\begin{equation}\label{RSVD}
{\mathbf X}\cong\overline{\mathbf U}{\mathbf S}{{\mathbf V}^T},
\end{equation}
where $\overline{\mathbf U}={\mathbf Q}{\mathbf U}$. Note that the size of ${\mathbf B}$ is much smaller than ${\mathbf X}$ and it is much easier to be handled.

Also, if we do not need the SVD of the original data matrix ${\mathbf X}$ and just a low-rank approximation of it is required, this can be easily achieved by the following approximation
\begin{equation}\label{QBAppro}
{\mathbf X}\cong{\mathbf Q}{\mathbf B},\,\,\,\,\,\,\,\,\,\,\,\,{\mathbf B}={\mathbf Q}^T{\mathbf X} \in \mathbb{R}^{R\times J}.
\end{equation}
Low-rank matrix approximation \eqref{QBAppro} is known as QB approximation of the matrix ${\mathbf X}$. 

\begin{rem}({\bf Oversampling and power iteration methods})\label{rempo}
The oversampling and power iteration techniques can be used to improve the solution accuracy of the randomized algorithms. In the oversampling technique, additional random vectors ($R+p$ random vectors instead of $R$ random vectors) are used to better capture the range of a matrix. Also, the power iteration technique is used when the decay rate of the singular values is not very fast. This technique replaces a matrix ${\mathbf X}$ by ${\left({\mathbf X}{{\mathbf X}^T}\right)^q}{\mathbf X}$ ($q$ is a nonnegative integer number) and randomized algorithms are applied to this new matrix. Considering the SVD, ${\mathbf X}={\mathbf U}{\mathbf S}{\mathbf V}^T$, we have ${\left({\mathbf X}\,{{\mathbf X}^T}\right)^q}\,{\mathbf X}={\mathbf U}\,{{\mathbf S}^{2q + 1}}\,{{\mathbf V}^T}$, and it is seen that the left and right singular vectors of the new matrix are the same as those of the matrix ${\mathbf X}$, but the singular values of the latter have a faster decay rate. This can improve the solution accuracy obtained by the randomized algorithms. Selections of the oversampling and the power iteration parameters depend on factors such as {\it the matrix dimensions}, {\it the singular spectrum decay} and {\it the random test matrices used in random projection}. However, if the Gaussian random matrices are used, it was experimentally confirmed that oversampling $p=5,10$ and power iteration $q=1,2$ are enough to achieve a good approximation \cite{halko2011finding}.
\end{rem}

The most computationally expensive and time-consuming part of the random projection technique is multiplying an original data matrix with random matrices. This costs $\mathcal{O}(IJR)$ flops. In order to accelerate this step, structured random matrices such as sparse random matrices \cite{achlioptas2003database, li2006very} and subsampled random Fourier transform (or SRFT) \cite{woolfe2008fast, liberty2009accelerated}, subsampled Hadamard transforms, sequence of Givens rotations \cite{liberty2009accelerated, ailon2009fast} can be used which are established techniques and used extensively in the literature \cite{rokhlin2008fast}. For example, the SRTF reduces the computation cost to $\mathcal{O}(IJ{\rm\,log}(R))$. These computations can be further reduced by exploiting compact random matrices \cite{sun2018tensor,batselier2018computing}. Randomized algorithms reduce the asymptotic complexity of deterministic algorithms for computation of the SVD from $\mathcal{O}\left( {IJR} \right)$ to $\mathcal{O}\left( {IJ\log (R) + (I + J){R^2}} \right)$ \cite{halko2011finding}.

The basic form of randomized SVD (BRSVD) algorithm equipped with the oversampling and the power iteration strategies is outlined in Algorithm \ref{ALG1} \cite{halko2011finding}. 

\RestyleAlgo{ruled}
\LinesNumbered
\begin{algorithm}
\SetKwInOut{Input}{Input}
\SetKwInOut{Output}{Output}  \Input{A data matrix ${\mathbf X}\in \mathbb{R}^{I\times J}$, target rank $R$, oversampling parameter $p$ and power iteration parameter $q$}  \Output{SVD factor matrices ${\mathbf U}\in\mathbb{R}^{I\times R},{\mathbf S}\in\mathbb{R}^{R\times R},{\mathbf V}\in\mathbb{R}^{R\times J}$}
	\caption{Basic Randomized SVD algorithm with oversampling and power iteration}\label{ALG1}
      Generate a random matrix ${\boldsymbol\Omega}\in\mathbb{R}^{J\times {(p+R)}}$ with prescribed probability distribution \\
      Form ${\mathbf Y}={{\left( {{\mathbf X}{{\mathbf X}^T}} \right)^q}{\mathbf X}}{\boldsymbol\Omega}$\\
      Compute QR decomposition ${\mathbf Y}={\mathbf Q}{\mathbf R}$\\
      Compute ${\mathbf B}={\mathbf Q}^T{\mathbf X}$\\
      Compute the full SVD of the matrix ${\mathbf B}=\overline{\mathbf U}\,\overline{{\mathbf S}}\,\overline{\mathbf V}^T$\\
      $\widetilde{\mathbf U}={\mathbf Q}\overline{\mathbf U}$\\
       ${\mathbf U}=\widetilde{\mathbf U}(:,1:R),\,{\mathbf S}=\overline{{\mathbf S}}(1:R,1:R),\,\,{\mathbf V}=\overline{\mathbf V}(:,1:R)$\\   	
\end{algorithm}
For stability issues, one should avoid computing expression $({\mathbf X}{\mathbf X}^T)^q{\mathbf X}$, explicitly and it should be computed sequentially using economic QR decomposition or LU decomposition \cite{halko2011finding}.

It has been proved that for oversampling parameter $P \ge 2$, power iteration $q$ and target rank $R$, the accuracy of approximation \eqref{RQR} is
\[
\begin{array}{l}
\mathbb{E}\left( \left\|{\cX - \cQ{\cQ^T}X}\right\|_2 \right) \le \\
\,\,\,\,\,\,\,\,\,\,{\left( {1 + \sqrt {\frac{R}{{P - 1}}}  + \frac{{e\sqrt {R + P} }}{P}\sqrt {\min (I,J) - k} } \right)^{1/(2q + 1)}}{\sigma _{k + 1}}
\end{array}
\]
where $\mathbb{E}$ denotes the mathematical expectation with the Gaussian random matrices and $\sigma_{R+1}$ is the $(R+1)$-th largest singular values of the matrix ${\mathbf X}$ \cite{martinsson2019randomized}. So, a nearly optimal low-rank approximation can be achieved using the randomization framework.


If two dimensions of an unfolding matrix are large, then it is possible to perform reduction on both sides. This can be performed by two random matrices. The structure of two-sided randomized algorithm is outlined in Algorithm \ref{TwoSide-SVD} \cite{halko2011finding}. Both Algorithms \ref{ALG1} and \ref{TwoSide-SVD} are randomized multi-pass algorithms because they compute multiplication with the original data matrix ${\mathbf X}$ in Line 4. These algorithms can be modified to become single-pass, e.g., by considering in Algorithm \ref{ALG1}
\begin{equation}\label{Onepass}
{\mathbf B} \cong {\left( {{\boldsymbol\Omega}_2}{\mathbf Q} \right)^\dag }{\mathbf W},\,\,\,\,\,{\mathbf W} = {{\boldsymbol\Omega}_2 {\mathbf X}},\,\,{\boldsymbol\Omega}_2\in\mathbb{R}^{R\times I},
\end{equation}
and in Algorithm \ref{TwoSide-SVD} 
\begin{equation}\label{Onepass2}
{\mathbf B} \cong {\left( {{{\boldsymbol{\Omega}}_2}\,{{\mathbf Q}^{(1)}}} \right)^\dag }\,{\mathbf W}\,{\left( {{{\mathbf Q}^{(2)\,T}}{\boldsymbol{\Omega}}_1} \right)^\dag },\,\,\,{\mathbf W} = {\boldsymbol{\Omega}}_2\,{\mathbf X}\,{{\boldsymbol{\Omega}}_1}.
\end{equation}
The benefit of these approaches is that they avoid computation of terms ${\mathbf Q}^T{\mathbf X}$ and ${\mathbf Q}_1^T{\mathbf X}{\mathbf Q}_2$ which may be computationally expensive, especially when the data matrix is stored out-of-cores where the cost of communication may exceed our main computations. Instead, in formulations \eqref{Onepass} and \eqref{Onepass2}, the original data matrix ${\mathbf X}$ is sketched by the random projection technique and the corresponding matrix ${\mathbf B}$ is obtained by solving some well-conditioned overdetermined linear least-squares problems \cite{tropp2017practical}. On the other hand, random matrix multiplication can be performed relatively fast by employing structured random matrices. 
We should note that this strategy passes the original data matrix ${\mathbf X}$ only once because all sketching procedures can be done in the first pass over the raw data. Other types of single-pass techniques can be found in \cite{tropp2017practical, halko2011finding, woodruff2014sketching, boutsidis2016optimal, upadhyay2018price, bajaj2019sketchycoresvd}. 
\RestyleAlgo{ruled}
\LinesNumbered
\begin{algorithm}
\SetKwInOut{Input}{Input}
\SetKwInOut{Output}{Output}  \Input{A data matrix ${\mathbf X}\in \mathbb{R}^{I\times J}$, target rank $R$}  \Output{SVD factor matrices ${\mathbf U}\in\mathbb{R}^{I\times R},{\mathbf S}\in\mathbb{R}^{R\times R},{\mathbf V}\in\mathbb{R}^{R\times J}$}
	\caption{Two-Sided Randomized SVD}\label{TwoSide-SVD}
	Draw prescribed random matrices ${{\boldsymbol\Omega}_1} \in {\mathbb{R}^{J \times R}},\,{{\boldsymbol\Omega}_2} \in {\mathbb{R}^{I \times R}}$\\
Compute ${{\mathbf Y}_1} = {\mathbf X}{{\boldsymbol\Omega}_1}$ and ${{\mathbf Y}_2} = {\mathbf X}^T{{\boldsymbol\Omega}_2}$\\
Compute economic QR decompositions ${{\mathbf Y}_1} = {{\mathbf Q}^{(1)}}{{\mathbf R}_1},\,{{\mathbf Y}_2} = {{\mathbf Q}^{(2)}}{{\mathbf R}_2}$\\
Compute ${\mathbf B} = {{\mathbf Q}^{(1)}}^T{\mathbf X}{{\mathbf Q}^{(2)}}$\\
Compute the SVD of matrix
${\mathbf B}={\overline{\mathbf U}}\,{\overline{\mathbf S}}\,{\overline{\mathbf V}^T}$\\
Set $\widetilde{\mathbf U}={\mathbf Q}^{(1)}{\overline{\mathbf U}}$ and $\widetilde{\mathbf V}={{\mathbf Q}^{(2)}\overline{\mathbf V}}$
\end{algorithm}

A main drawback of the random projection algorithms described so far is that they need an estimation of the matrix rank which may be a difficult assumption. To resolve this difficulty, randomized {\em rank-revealing algorithms}, or equivalently {\em randomized fixed-precision algorithms} have been developed in \cite{martinsson2016randomizedReveal, yu2018efficient}, where for a given approximation error bound, the rank and corresponding low-rank approximation can be computed automatically.

\subsection{Sampling Techniques}\label{Sec:RS}
The other type of randomized algorithms for finding low-rank approximations of matrices is based on selecting a part of its columns or rows randomly. 
The previous approach does not choose individual columns of the matrix $\mathbf X$, but instead it finds a set of linear combinations of columns (rows) capturing the column (row) space of the matrix $\mathbf X$ with high probability. Sampling technique is an alternative approach to the random projection where instead of multiplying the original data matrix with random matrices, some columns or rows of it are selected and the original data matrix is compressed in this manner\footnote{Component selection is another sampling approach where some components (not necessarily fibers) of the original data matrix are selected using some probability distributions \cite{drineas2016randnla}. This is also known as sparsification procedure.}. The procedure of column or
row selection can be performed based on different kinds of probability distributions such as uniform distribution, leverage scores \cite{drineas2008relative}, length squared \cite{frieze2004fast, Frieze1998} and also with/without replacement\footnote{It has been reported that for a fixed distribution, and when the columns are uniformly correlated, the uniform sampling without replacement works quite well in practice \cite{drinea2001randomized, kumar2009sampling,kumar2012sampling}.}. The accuracies of results highly depend on the probability distribution used in the sampling procedure and either {\em additive} or {\em relative} approximation error can be achieved. 

Let ${\mathbf X}\in \mathbb{R}^{I\times J}$ be a given data matrix, and ${\mathbf C}$ be a matrix containing the selected columns, then we have two kinds of approximation as\footnote{The same concepts are used for spectral norm ${\left\| . \right\|_2}$.} \cite{mahoney2011randomized}
\begin{itemize}
\item {\bf Relative-error approximation}\footnote{It is also called multiplicative error.}
\begin{equation}\label{RelErr}
\begin{array}{l}
{\left\| {{\mathbf X} - {\mathbf C}{{\mathbf C}^\dag}{\mathbf X}} \right\|_F} \le (1+\epsilon){\left\| {{\mathbf X} - {{\mathbf X}_R}} \right\|_F},
\end{array}
\end{equation}

\item {\bf Additive-error approximation} 
\[
\begin{array}{l}
{\left\| {{\mathbf X} - {\mathbf C}{{\mathbf C}^\dag}{\mathbf X}} \right\|^2_F} \le {\left\| {{\mathbf X} - {{\mathbf X}_R}} \right\|^2_F} + \epsilon\left\| {{{\mathbf X}_R}} \right\|^2_F.
\end{array}
\]
\end{itemize}
where ${\mathbf X}_R$ is the best rank $R$-approximation of the matrix $\mathbf X$.
Clearly, randomized algorithms with relative-error guarantees are of much more interest than the additive-error ones. The number of selected columns should be large enough to capture the range of a matrix and satisfy 
the relation \eqref{RelErr}.

The best existing random sampling algorithms use the leverage score probability distribution \cite{drineas2012fast} while their computation is expensive because of the computation of the SVD. In \cite{drineas2012fast}, a fast and computationally efficient approach has been proposed for computing the leverage scores. In some applications, we need interpretable low-rank matrix approximations, and columns should be taken from the original data matrix. Here, in the first step, the leverage scores of an underlying data matrix can be computed by the randomized SVD and then the columns of the original data matrix are selected based on the leverage score probability distribution. Note that exploiting random uniform sampling of columns (fibers) provides relatively fast computation. However, this method works well only if a data matrix has low coherence \cite{mahoney2011randomized} otherwise they do not provide a reliable compression. An alternative approach is to spread out the information of the data matrix or uniformize it by a prior random projection. This preprocessing step allows applying the uniform sampling\cite{drineas2006fast}. 

The sampling low-rank matrix approximation is outlined in Algorithm \ref{ColSamp2}. It picks $R$ columns of the matrix $\mathbf X$ and constructs the matrix $\mathbf Y$ \cite{drineas2006fast}, where in Line 3 of Algorithm \ref{ColSamp2}, the selected columns are scaled \cite{drineas2006fast}. This scaling procedure provides an unbiased estimator for the matrix-matrix multiplication \cite{mahoney2011randomized}.

\RestyleAlgo{ruled}
\LinesNumbered
\begin{algorithm}
      \SetKwInOut{Input}{Input}
\SetKwInOut{Output}{Output}  \Input{A data matrix ${\mathbf X} \in {\mathbb{R}^{I \times J}},\,R \in {\mathbb{Z}^+}$ such that $1\le R \le J,\,\left\{ {{p_j}} \right\}_{j = 1}^J,\,$ ${\sum\limits_{j = 1}^J {{p_j}}=1}\,\,{\rm{and}}$ ${p_j} \ge 0$}  \Output{Low-rank matrix approximation ${\mathbf X} \cong {\mathbf Q}\,{\mathbf B},\,{\mathbf Q}\in\mathbb{R}^{I\times R}\,\,{\mathbf B}\in\mathbb{R}^{R\times J}$}
	\caption{Sampling algorithm for low-rank matrix approximation}\label{ColSamp2}
       \For{$r=1,2,\ldots,R$}
      {
      Pick$\,\,\,{j_r} \in \{1, 2, \ldots ,J\}\,\,\,$with$\,\,\Pr \left( {{j_r} = \alpha } \right) = {p_\alpha }$ $\alpha  = 1, 2, \ldots ,J$\\
      Set$\,\,\widetilde{\mathbf Y}(:,j_r) = \frac{{\mathbf X}(:,j_r)}{\sqrt {R{p_{{j_r}}}}}$\\
 	}
 	Compute ${\widetilde{\mathbf Y}}={\mathbf Q}\,{\mathbf R},\,\,\,{{\mathbf B}} ={\mathbf Q}^T{\mathbf X}$
\end{algorithm}

\subsection{count-sketch}\label{count-sketch}
The random projection approach needs expensive multiplication of an original data matrix with random matrices. 
Count-sketch uses an alternative trick to generate a linear combinations of the columns without multiplication with random matrices. 
Count-sketch was originally introduced in the context of data stream \cite{charikar2002finding} and has been used 
to speed-up the matrix-matrix computations \cite{pham2013fast, pagh2013compressed, clarkson2017low}. Similar to the  random projection and the sampling approaches, the {\em count-sketch} technique aims at capturing the range of a given matrix. It consists of the following three steps:
\begin{itemize}
\item Hashing procedure,

\item Grouping together the columns with the same hash numbers,

\item Signing the columns and summing each group as a representative column.

\end{itemize}
Assume that we want to compute a rank-$R$ approximation of a matrix. In the hashing procedure, all columns of the data matrix are labeled with numbers $1,2,\ldots,R$ uniformly, i.e., with probability $1/R$.

In the second step, the columns with the same label are grouped together and after signing the columns\footnote{By signing we mean multiplying the columns by $\pm 1$ uniformly. This is also called binary Rademacher variable.}, the columns in each group are summed up as representative columns. 

Clearly, the cost of the first step is negligible, and the main operation in the second step is summing the columns (with integer coefficients $\pm 1$). Compared to multiplying the original data matrix by a random matrix, this technique can significantly reduce the running time. 
\subsection{Randomized Algorithms for Solving Least-Squares Problems}
In this section, we briefly describe the randomized algorithms for solving large-scale least-squares problems. 
In the procedure of computation of the Canonical Polyadic Decomposition (CPD) \cite{hitchcock1927expression, hitchcock1928multiple} and the Tucker decomposition, we often need to solve least-squares problems with tall and skinny coefficient matrices, i.e., the number of rows is much more than the number of columns. For example, the following two matrices 
\begin{eqnarray}\label{KHM}
\mathop  {\mathlarger\odot} \limits_{i = 1}^N {{\mathbf Q}^{\left( n \right)}}\in{\mathbb{R}^{\left( {\prod\limits_{n = 1}^N {{I_n}} } \right) \times {R_n}}},
\end{eqnarray}
\begin{equation}\label{KM}
\mathop  {\mathlarger\otimes} \limits_{n = 1}^N {{\mathbf Q}^{\left( n \right)}}\in{\mathbb{R}^{\left( {\prod\limits_{n = 1}^N {{I_n}} } \right) \times \left( {\prod\limits_{n = 1}^N {{R_n}} } \right)}},
\end{equation}
have tall and skinny structures where ${{\mathbf Q}^{\left( n \right)}} \in {\mathbb{R}^{{I_n} \times {R_n}}},\,\,\,{I_n} \gg {R_n}$, 
and arise in the procedure of computation of the CPD and the Tucker decomposition \cite{kolda2009tensor}.

\begin{figure}
\begin{center}
\includegraphics[width=7.5cm,height=4cm]{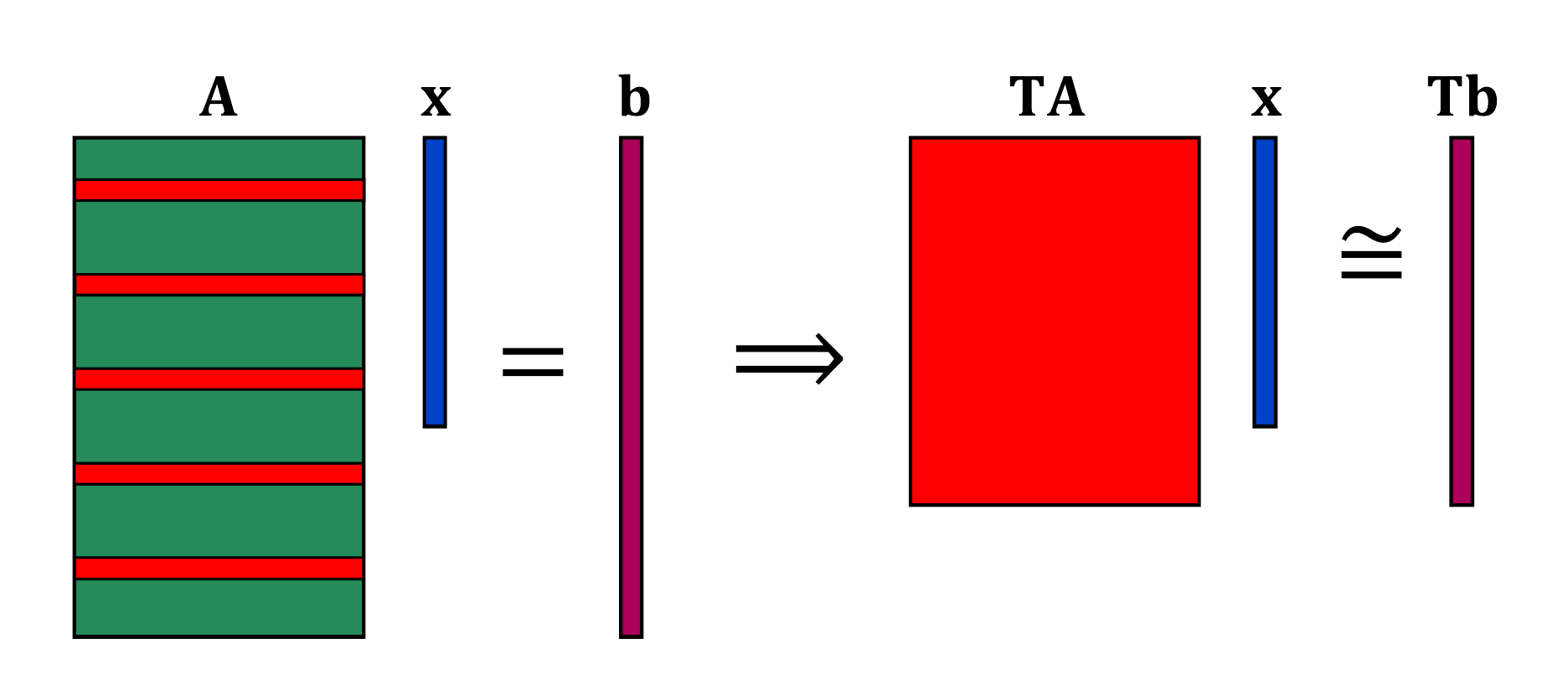}\\
\caption{Illustration of randomized row sampling technique for solving an overdetermined least-squares problem.} 
\label{FIGTTTR}
\end{center}
\end{figure}

Consider the following least-squares problem 
\begin{equation}\label{LST}
{\mathbf x} = \mathop {\arg \min }\limits_{\mathbf x} \,\,\left\| {{\mathbf A}{\mathbf x} - {\mathbf b}} \right\|_2,
\end{equation}
where the matrix ${\mathbf A}$ is either \eqref{KHM} or \eqref{KM}, the right-hand side is ${\mathbf b}\in\mathbb{R}^{\prod\limits_{n = 1}^N {{I_n}}}$ and the vector ${\mathbf x}\in\mathbb{R}^{R_n}$ (or ${\mathbf x}\in\mathbb{R}^{\prod\limits_{n = 1}^N {{R_n}}}$) is the unknown vector needs to be determined. The randomized or sketching techniques replace the least-squares problem \eqref{LST} by
\begin{equation}\label{RLST}
\tilde {\mathbf x} = \mathop {\arg \min }\limits_{\mathbf x} \,\,\left\| {{\mathbf T}{\mathbf A\mathbf x} - {\mathbf T}{\mathbf b}} \right\|_2,
\end{equation}
where \[
{\mathbf T}:{\mathbb{R}^{\prod\limits_{n = 1}^N {{I_n}} }} \to {\mathbb{R}^L},\,\,\,L \le \prod\limits_{n = 1}^N {{I_n}},
\] 
is a map which reduces the dimensionality of the matrix ${\mathbf A}$ and the right-hand side ${\mathbf b}$, see Figure \ref{FIGTTTR} for a graphical illustration. 
The transformation ${\mathbf T}$ can be the random projection, the sampling or the count-sketch technique discussed in Section \ref{Sec:RandomMatrix}. 
It is expected to find a solution $\tilde {\mathbf x}$ of the least-squares problem \eqref{RLST}  which satisfies  
\[
{\left\| {{\mathbf A}\tilde {\mathbf x} - {\mathbf b}} \right\|_2} \le \left( {1 + \varepsilon } \right){\left\| {{\mathbf Ax^{*}} - {\mathbf b}} \right\|_2},
\]
with a high probability where ${\mathbf x}^*$ is the optimal solution of least-squares problem \eqref{LST} and $\epsilon>0$ is a given tolerance \cite{woodruff2014sketching}. 

In next section, we introduce the Tucker decomposition and the HOSVD along with algorithms for their computations.

\section{Tucker Decomposition and Higher Order Singular Value Decomposition (HOSVD)}\label{Sec:HOSVD}
Let $\underline{\mathbf X}\in{\mathbb{R}^{{I_1} \times I_2\times  \cdots  \times {I_N}}}$, then the Tucker decomposition of the tensor $\underline{\mathbf X}$ admits the following model \cite{tucker1963implications, tucker1964extension, tucker1966some} 
\begin{equation}\label{HOSVD}
\underline{\mathbf X} \cong \underline{\mathbf S}{ \times _1}{{\mathbf Q}^{\left( 1 \right)}}{ \times _2}{{\mathbf Q}^{\left( 2 \right)}} \cdots { \times _N}{{\mathbf Q}^{\left( N \right)}},
\end{equation}
where $\underline{\mathbf S} \in {\mathbb{R}^{{R_1} \times  R_2\times \cdots  \times {R_N}}}$ is a core tensor and ${{\mathbf Q}^{\left( n \right)}} \in {\mathbb{R}^{{I_n} \times {R_n}}},\,{R_n} \le {I_n},\,\,n=1,2,\ldots,N$ are factor matrices. A shorthand notation for the Tucker decomposition is
\[
\underline{\mathbf X} \cong \left[\kern-0.15em\left[ {\underline{\mathbf S};{{\mathbf Q}^{(1)}},{{\mathbf Q}^{(2)}}, \ldots ,{{\mathbf Q}^{(N)}}}
 \right]\kern-0.15em\right].
\]

\begin{figure}
\begin{center}
\includegraphics[width=8 cm,height=4 cm]{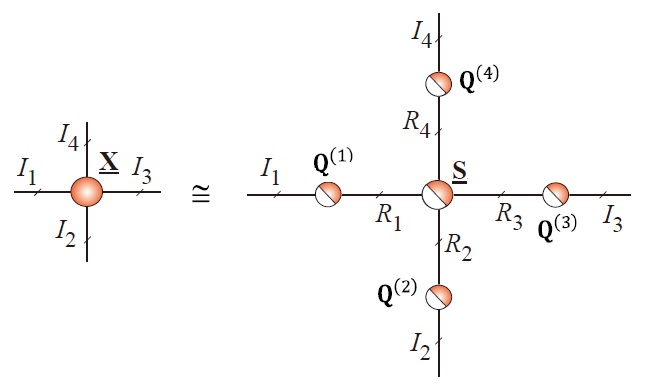}\\
\includegraphics[width=8 cm,height=4 cm]{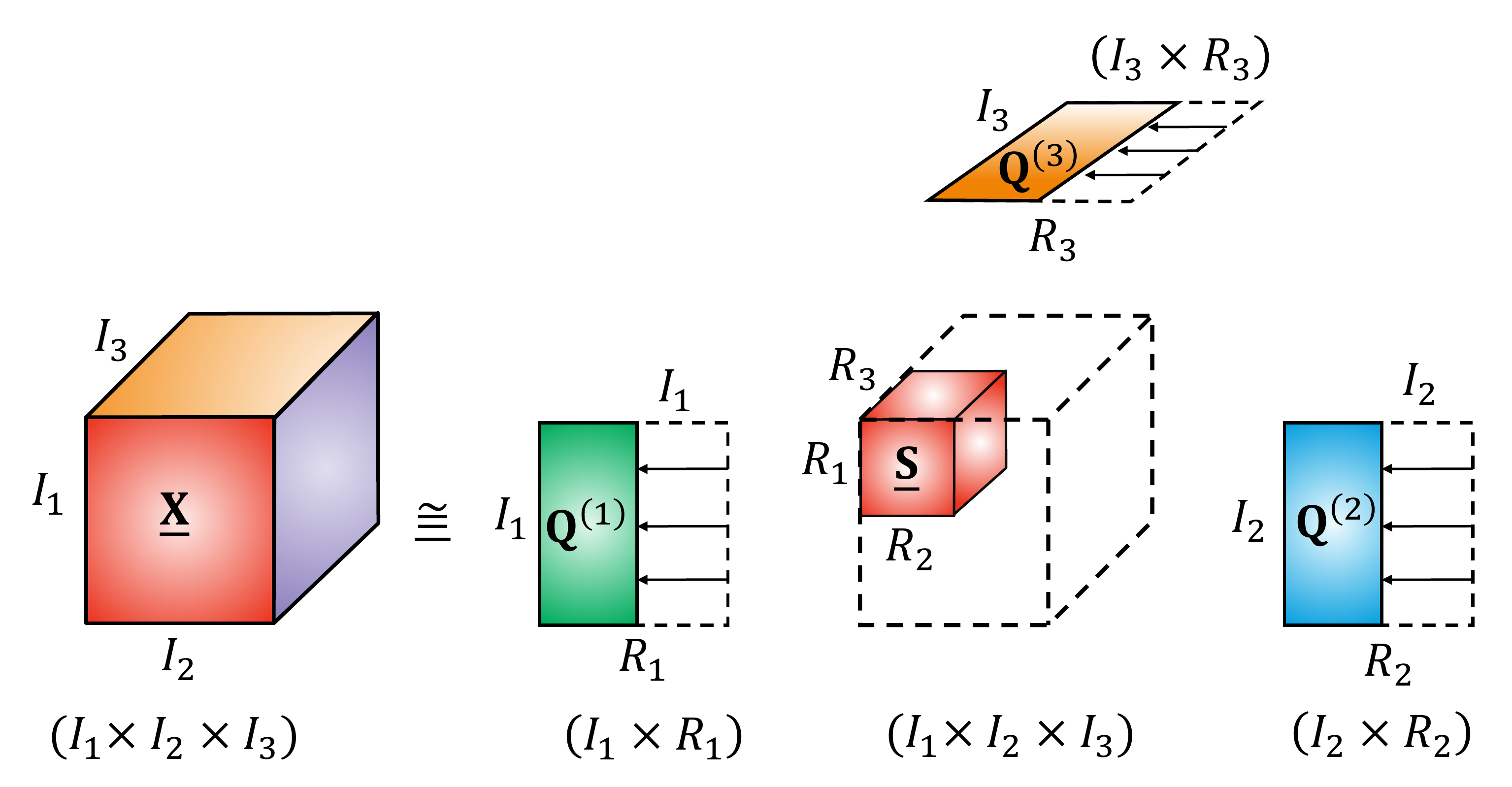}\\
\caption{\small{Graphical illustrations of the truncated HOSVD, upper is for a fourth order tensors and bottom is for a third order tensors.}}\label{HOSVDGraph}
\end{center}
\end{figure}

HOSVD or equivalently Multilinear SVD (MLSVD) is a constrained Tucker decomposition \cite{de2000multilinear} that ensures the orthogonality of factor matrices and all-orthogonality of the core tensor\footnote{A tensor is called all-orthogonal if all its slices in each mode are mutually orthogonal \cite{de2000multilinear}.}. For a graphical illustration of the HOSVD for a 3rd and a 4th-order tensors, see Figure \ref{HOSVDGraph}. 

Unlike the SVD, the core tensor $\underline{\mathbf S}$, in general, is not diagonal or even not nonnegative, but it has {\em pseudo-diagonal} property, which means that the Frobenius norm of slices in each mode is non-increasing as the index is increased \cite{de2000multilinear}. Also, another interpretation of this concept is  that the intensity (absolute value) of the components of the core tensor $\underline{\mathbf S}$ is mainly concentrated on one of the corners (position $(1,1,\ldots,1)$) and it decreases as the components move away from this corner. 


Due to the pseudo-diagonality property of the HOSVD which plays a role similar to the singular values in the SVD, a truncated version of the HOSVD is used in practice. To this end, first the orthogonal matrices ${\mathbf Q}^{(n)}\in\mathbb{R}^{{I_n} \times {R_n}},\,\,n=1,2,\ldots,N$ in \eqref{HOSVD} are computed through the truncated SVD (tSVD) or the randomized SVD (rSVD) algorithms for low-rank approximation of the unfolding matrices as
\begin{equation}\label{UnfoldHOSVD}
{{\mathbf X}_{\left( n \right)}} \cong {{\mathbf Q}^{\left( n \right)}}{{\mathbf S}_n}{{\mathbf V}^{\left( n \right)\,\,T}} \in {\mathbb{R}^{{I_n} \times {I_1} \cdots {I_{n - 1}}{I_{n + 1}} \cdots {I_N}}},
\end{equation}
where ${\mathbf V}^{(n)}\in\mathbb{R}^{{I_1\cdots I_{n-1}I_{n+1}\cdots I_N}\times R_n},\,{\mathbf S}_n\in\mathbb{R}^{R_n\times R_n}$ and $R_n$ is the numerical rank of the unfolding matrix ${\mathbf X}_{(n)}$. Then, due to the orthogonality of the factor matrices ${\mathbf Q}^{(n)}\in\mathbb{R}^{{R_n} \times {I_n}},\,\,n=1,2,\ldots,N$, the core tensor can be computed as
\begin{equation}\label{Core}
\underline{\mathbf S} = \underline{\mathbf X}{ \times _1}{{\mathbf Q}^{\left( 1 \right)\,T}}{ \times _2}{{\mathbf Q}^{\left( 2 \right)\,T}} \cdots { \times _N}{{\mathbf Q}^{\left( N \right)\,T}}\in\mathbb{R}^{R_1\times R_2\times \cdots\times R_N}.
\end{equation}

Since, we only need the left singular vectors in \eqref{UnfoldHOSVD}, the eigenvalue decomposition (EVD) of the Gramian matrix ${{\mathbf X}_{(n)}}{\mathbf X}_{(n)}^T = {\mathbf Q}^{(n)}{\mathbf\Sigma}^2{{\mathbf Q}^{(n)\,T}}$ can also be used in our computations. The Crank-Nicholson-like algorithm \cite{phan2014fast,camarrone2018fast}, rank revealing QR (LU) decompositions \cite{chan1987rank, pan2000existence} or Alternating Least Squares (ALS) type techniques \cite{xiao2020efficient, li2020tucker} are other possible technique for computation of the factor matrices. 

The $N$-tuple $(R_1,R_2,\ldots,R_N)$ is called multilinear rank or Tucker rank of tensor $\underline{\mathbf X}$.

The truncated HOSVD does not provide the best multilinear rank approximation in the least-squares sense while a quasi-best approximation can be achieved \cite{de2000multilinear} as follows
\[
\left\| {\underline{\mathbf X} - \underline{\mathbf S}{ \times _1}{{\mathbf Q}^{\left( 1 \right)}}{ \times _2} \cdots { \times _N}{{\mathbf Q}^{\left( N \right)}}} \right\|_F \le \sqrt N \left\| {\underline{\mathbf X} - {{\underline{\mathbf X}_{Best}}}} \right\|_F,
\]
where $\underline{\mathbf X}_{best}$ is the best multilinear rank approximation of the tensor $\underline{\mathbf X}$.

Substituting \eqref{Core} in \eqref{HOSVD}, we have
\begin{equation}\label{TK1}
\underline{\mathbf X} \cong \underline{\mathbf X}{ \times _1}{{\mathbf Q}^{\left( 1 \right)}}{{\mathbf Q}^{\left( 1 \right)\,T}}{ \times _2}{{\mathbf Q}^{\left( 2 \right)}}{{\mathbf Q}^{\left( 2 \right)\,T}} \cdots { \times _N}{{\mathbf Q}^{\left( N \right)}}{{\mathbf Q}^{\left( N \right)\,T}}.
\end{equation}
Note that if the factor matrices ${\mathbf Q}^{(n)},\,n=1,2,\ldots,N$ are not orthogonal and just have full column rank, then \eqref{Core} and \eqref{TK1} are replaced by\footnote{If $\mathbf Q$ is full-rank then ${{\mathbf Q}^\dag }{\mathbf Q} = {\mathbf I}$, where $\mathbf I$ is the identity matrix of conforming dimension.}  
\begin{equation}\label{Core2}
\underline{\mathbf S}\cong \underline{\mathbf X}{ \times _1}{{\mathbf Q}^{\left( 1 \right)\dag}}{ \times _2}{{\mathbf Q}^{\left( 2 \right)\dag}} \cdots { \times _N}{{\mathbf Q}^{\left( N \right)\dag}},
\end{equation}
and 
\begin{equation}\label{AP}
\underline{\mathbf X}\cong \underline{\mathbf X}{ \times _1}{{\mathbf Q}^{\left( 1 \right)}}{{\mathbf Q}^{\left( 1 \right)\,\dag}}{ \times _2}{{\mathbf Q}^{\left( 2 \right)}}{{\mathbf Q}^{\left( 2 \right)\,\dag}} \cdots { \times _N}{{\mathbf Q}^{\left( N \right)}}{{\mathbf Q}^{\left( N \right)\,\dag}},
\end{equation}
where $\dag$ is the Moore-Penrose pseudoinverse operator.
It is worth mentioning that if the dimension of a particular mode, say mode $n$, is relatively small, then we can ignore reduction in the mentioned mode or, equivalently, ignore its corresponding orthogonal matrix ${\mathbf Q}^{(n)}$ and remove it from \eqref{TK1} and \eqref{AP}.
In view of \eqref{TK1} and \eqref{AP}, it is seen that to find an approximate HOSVD or an approximate Tucker decomposition, good approximations for the range of unfolding matrices are required. More precisely, the main problem is how to find a good approximation for the range of unfolding matrices ${\mathbf X}_{(n)}$, i.e., ${\mathbf Q}^{(n)}$. This problem can be formally formulated as follows:

Given a data tensor $\underline{\mathbf X}\in{\mathbb{R}^{{I_1} \times  I_2\times \cdots  \times {I_N}}}$ and tolerance $\epsilon$, the problem is finding orthogonal matrices ${\mathbf Q}^{(n)}\in\mathbb{R}^{{I_n} \times {R_n}}$, satisfying
\begin{equation}\label{DT1}
{\left\| {\underline{\mathbf X} - \underline{\mathbf X}{ \times _1}{{\mathbf Q}^{\left( 1 \right)}}{{\mathbf Q}^{\left( 1 \right)\,\,T}} \cdots { \times _N}{{\mathbf Q}^{\left( N \right)}}{{\mathbf Q}^{\left( N \right)\,\,T}}} \right\|_F} \le {\sqrt N}\varepsilon,
\end{equation}
and in the case when ${\mathbf Q}^{(n)}$ are not orthogonal
\begin{equation}\label{DT2}
{\left\| {\underline{\mathbf X} - \underline{\mathbf X}{ \times _1}{{\mathbf Q}^{\left( 1 \right)}}{{\mathbf Q}^{\left( 1 \right)\dag}} \cdots { \times _N}{{\mathbf Q}^{\left( N \right)}}{{\mathbf Q}^{\left( N \right)\dag}}} \right\|_F} \le {\sqrt N} \varepsilon.
\end{equation}
We will return to these problems later in Section \ref{Sec:RHOSVD}.
\subsection{Sequentially Truncated HOSVD (STHOSVD) Algorithm}
In formulation \eqref{UnfoldHOSVD}, low-rank approximations of the unfolding matrices ${\mathbf X}_{(n)}$ are required where all unfolding matrices have the same number of elements as the original data tensor $\underline{\mathbf X}$. A more efficient algorithm with a less computational complexity for computing the HOSVD is {\em Sequentially Truncated HOSVD} (STHOSVD) algorithm \cite{vannieuwenhoven2012new}. Instead of dealing with all unfolding matrices (all with the same number of elements as the original data tensor $\underline{\mathbf X}$) and finding their low-rank matrix approximations, this algorithm reduces the size of underlying unfolding matrices sequentially. More precisely, at each iteration of the algorithm, a new unfolding matrix which approximates a specific factor matrix has a smaller size than the previous ones. 
This results in significant speed-up with comparable accuracy or sometimes even better accuracy\footnote{Although there is no theoretical evidence establishing better accuracy of SHOSVD compared to the HOSVD, numerical simulations highly support this. However, there are some  examples for which the HOSVD provides better accuracy compared to the SHOSVD \cite{vannieuwenhoven2012new}.}. This procedure is described in Algorithm \ref{SHOSVD} \cite{vannieuwenhoven2012new}. It is worth mentioning that in Lines 4 of Algorithm \ref{SHOSVD}, an $n$-mode (matrix) product of the tensor $\underline{\mathbf S}$ by the matrix ${\mathbf Q}^{{(n)}\,T}$ is required. This can be performed in an equivalent way. To be more precise, let $\left[ {{{\mathbf Q}^{\left( n \right)}},{{\mathbf \Lambda}^{\left( n \right)}},{{\mathbf V}^{\left( n \right)}}} \right]$ be the factor matrices of the truncated SVD of the unfolding matrix ${\mathbf S}_{(n)}$. If $\underline{\mathbf Y} \leftarrow \underline{\mathbf S}{ \times _n}{{\mathbf Q}^{\left( n \right)\,\,T}}$ then it is easy to show that ${{\mathbf Y}_{\left( n \right)}} = {{\mathbf \Lambda}^{\left( n \right)}}\,{{\mathbf V}^{\left( n \right)\,\,T}}$. This means that in Algorithm \ref{SHOSVD}, in Line 3, if we compute all factor matrices $\left[ {{{\mathbf Q}^{\left( n \right)}},{{\mathbf \Lambda}^{\left( n \right)}},{{\mathbf V}^{\left( n \right)}}} \right]$, then Line 4 can be replaced by ${{\mathbf S}_{\left( n \right)}} = {{\mathbf \Lambda}^{\left( n \right)}}\,{{\mathbf V}^{\left( n \right)\,\,T}}$.

%

Another advantage of the STHSOVD compared to the THOSVD is that in the STHOSVD algorithm, the core tensor is automatically computed by the algorithm in the last step while in the THOSVD algorithm, all factor matrices are first computed and then the core tensor is computed through formulation \eqref{Core} or \eqref{Core2}. This may lead to the intermediate data explosion phenomenon \cite{kolda2008scalable}. Note that Algorithm \ref{SHOSVD}, starts from the first factor matrix and other factor matrices are computed in ascending order, i.e. $p=[1,2,\ldots,N]$, but other orderings are also possible. It turns out that the ordering used within the STHOSVD algorithm affects the approximation accuracy and speed-up \cite{vannieuwenhoven2012new}.

\RestyleAlgo{ruled}
\begin{algorithm}
\LinesNumbered
\SetKwInOut{Input}{Input}
\SetKwInOut{Output}{Output}  \Input{A data tensor $\underline{\mathbf X} \in {\mathbb{R}^{{I_1}\times I_2\times \cdots  \times {I_N}}}$ and a tolerance $\epsilon$;}  \Output{Approximative HOSVD of the tensor $\underline{\mathbf X}$ $\underline{\mathbf X} \cong \left[\kern-0.15em\left[ {\underline{\mathbf S};{{\mathbf Q}^{(1)}},{{\mathbf Q}^{(2)}}, \ldots ,{{\mathbf Q}^{(N)}}}
 \right]\kern-0.15em\right]$ and multilinear rank $(R_1,R_2,\ldots,R_N)$}
	\caption{Sequentially Truncated HOSVD (STHOSVD) Algorithm} \label{SHOSVD}
Set $\underline{\mathbf S}=\underline{\mathbf X}$\\
\For{$n=1,2,\ldots,N$}
{$[{\mathbf Q}^{(n)},\mathtt{\sim},\mathtt{\sim}]$ =  truncated-svd$\left({\mathbf S}_{\left(n\right)},\frac{\epsilon}{N}\right)$\\
$\underline{\mathbf S}=\underline{\mathbf S}\times_{n} {\mathbf Q}^{(n)\,T}$\\
}
\end{algorithm}

\subsection{Higher Order Orthogonal Iteration (HOOI) Algorithm}

Unlike {\em Eckart-Young} property for matrices \cite{eckart1936approximation}, neither the truncated HOSVD (THOSVD) nor the sequentially truncated HOSVD (STHOSVD) provides the best multilinear rank approximation for higher order tensors. Several algorithms have been 
developed for the computation of the best multilinear rank approximation among which we can mention Newton-Grassmann algorithm \cite{elden2009newton}, Riemannian trust-region algorithm \cite{ishteva2011best}, Higher-Order Orthogonal Iteration (HOOI) algorithm\footnote{This algorithm is also known as Tucker-ALS algorithm.} \cite{de2000best} etc. HOOI algorithm is the simplest kind of these algorithms which is based on the idea of ALS technique where at each iteration of the algorithm, all factor matrices are fixed except one and it is updated by solving a least-squares problem. To be more precise, consider the Tucker decomposition \eqref{HOSVD}
and apply the $n$-unfolding on both sides of \eqref{HOSVD}, as
\[
\begin{array}{l}
{{\mathbf X}_{\left( n \right)}} = {{\mathbf Q}^{\left( n \right)}}{{\mathbf S}_{\left( n \right)}}\left( {{{\mathbf Q}^{\left( N \right)}} \otimes\cdots\otimes{{\mathbf Q}^{\left( {n + 1} \right)}}{\otimes}{{\mathbf Q}^{\left( {n - 1} \right)}}} {\otimes}\right.\\
\hspace{6.3cm}\left. { \cdots {\otimes }{{\mathbf Q}^{\left( 1 \right)}}} \right)^T,
\end{array}
\]
for $n=1,2,\ldots,N$. Assume that the core tensor $\underline{\mathbf S}$ and all but one factor matrices $\left\{ {{{\mathbf Q}^{\left( p \right)}}} \right\},\,\,p \ne n$, are known. The unknown factor matrix ${\mathbf Q}^{(n)}$ can be computed by solving the following least-squares problem\footnote{In \cite{chachlakis2019l1}, the $l_1$-norm form of least-squares problem \eqref{LSS} is solved.}
\begin{eqnarray}\label{LSS}
{{\mathbf Q}^{\left( n \right)}} = \mathop {\arg \min }\limits_{{\mathbf Q} \in {\mathbb{R}^{{I_n} \times {R_n}}}} \,\,\left\| {{\mathbf A}^{(n)}{{\mathbf Q}^T} - {\mathbf X}_{\left( n \right)}^T} \right\|_F^2,
\end{eqnarray}
where 
\[
{\mathbf A}^{(n)}=\left( {\mathop  {\mathlarger\otimes} \limits_{p = N,p \ne n}^1 {{\mathbf Q}^{\left( p \right)}}} \right){\mathbf S}_{\left( n \right)}^T.
\]
It can be shown that the solution to the least-squares problem \eqref{LSS} is equivalent to finding $R_n$ leading left singular vectors of ${\mathbf Z}_{(n)}$ \cite{de2000multilinear} where 
\begin{equation}\label{EXpHO}
\underline{\mathbf Z} = \underline{\mathbf X} { \times _{p \ne n}}\left\{ {{{\mathbf Q}^{\left( p \right)}}} \right\}^T.
\end{equation}
This procedure is repeated for all factor matrices ${\mathbf Q}^{(n)}\in\mathbb{R}^{I_n\times R_n},\,n=1,2,\ldots,N,$ sequentially and in an usual Gauss-Seidel manner. After computing all factor matrices, the core tensor $\underline{\mathbf S}$ is updated in the following form
\begin{equation}\label{CUpdate}
\underline{\mathbf S} = \underline{\mathbf X}{ \times _1}{{\mathbf Q}^{\left( 1 \right)\,T}} \cdots { \times _N}{{\mathbf Q}^{\left( N \right)\,T}}.
\end{equation}
Also the computation of the core tensor $\underline{\mathbf S}$ in \eqref{CUpdate}, can be equivalently computed by solving the following least-squares problem 
\begin{equation}\label{LSE2}
\underline{\mathbf S} = \mathop {\arg \min }\limits_{ \underline{\mathbf S} \in {\mathbb{R}^{{R_1} \times {R_2} \times \cdots  \times {R_N}}}} \,\,\left\| {{\mathbf B}\,{{\mathbf s}{\left( : \right)}} - {{\mathbf x}{\left( : \right)}}} \right\|_2^2,
\end{equation}
where 
\[
{\mathbf B}={\mathop\otimes \limits_{n = N}^1 {{\mathbf Q}^{\left( n \right)}}},
\]
and ${\mathbf x}(:)$ and ${\mathbf s}(:)$ are the vectorization of the tensors $\underline{\mathbf X}$ and $\underline{\mathbf S}$, respectively.
This procedure is outlined in Algorithm \ref{AHOOI} \cite{de2000best}. 

\begin{center}
\RestyleAlgo{ruled}
\begin{algorithm}
\LinesNumbered
\SetKwInOut{Input}{Input}
\SetKwInOut{Output}{Output}
	  \Input{A data tensor $\underline{\mathbf X} \in {\mathbb{R}^{{I_1}\times I_2 \times \cdots  \times {I_N}}}$, and a multilinear rank $\left( {{R_1}, R_2,\ldots ,{R_N}} \right)$}
	        \Output{Approximative Tucker representation of the tensor $\underline{\mathbf X}$ as $\underline{\mathbf X} \cong \left[\kern-0.15em\left[ {\underline{\mathbf S};{{\mathbf Q}^{(1)}},{{\mathbf Q}^{(2)}}, \ldots ,{{\mathbf Q}^{(N)}}}
 \right]\kern-0.15em\right]$}
	        Initialize factor matrices ${\mathbf Q}^{(2)},\ldots,{\mathbf Q}^{(N)}$ with the Truncated HOSVD or random matrices\\
	        \While{A stopping criterion is not satisfied}{
  \For{$n = 1,2,\ldots, N$}
  {  $\underline{\mathbf Z} = \underline{\mathbf X}{ \times _{p \ne n}}\left\{ {{{\mathbf Q}^{\left( p \right)\,\,T}}} \right\}$;\\
${\mathbf Q}^{(n)}$ is constructed from $R_n$ leading left singular vectors of ${\mathbf Z}_{(n)}$
  }
  $\underline{\mathbf S} = \underline{\mathbf Z}{ \times _1}{{\mathbf Q}^{\left( 1 \right)\,\,T}}{ \times _2}{{\mathbf Q}^{\left( 2 \right)\,\,T}}\cdots {\times _N}{{\mathbf Q}^{\left( N \right)\,\,T}}$
  }
 \caption{Higher-Order Orthogonal Iteration (HOOI) algorithm}\label{AHOOI}
\end{algorithm}
\end{center}

\section{Randomized Tucker decomposition and randomized HOSVD}\label{Sec:RHOSVD}
It can be seen that the most computationally expensive part of the HOSVD or the STHOSVD is low-rank matrix approximation of the unfolding matrices. In randomized Tucker algorithms or randomized HOSVD, the truncated (or economic) SVD is replaced by the randomized low-rank matrix algorithms. In this section, we discuss a variety of randomized algorithms for computation of the Tucker decomposition and the HOSVD. 
\subsection{\bf Randomized random projection Tucker Decomposition}
Exploiting the idea of random projection, the randomized variants of the HOSVD and the HOOI algorithms can be straightforwardly developed. These algorithms are outlined in Algorithms \ref{RHOSVD} and \ref{HOOIRAND} \cite{zhou2014decomposition, wolf2019low}.  

In Algorithms \ref{RHOSVD}, the unfolding matrix is not necessary to be  computed explicitly. In \cite{che2021randomized}, instead of ${{\mathbf Z}_{\left( n \right)}}{{\mathbf \Omega} ^{\left( n \right)}},$ it is suggested to perform ${\underline {\mathbf W}}={\underline{\mathbf Z}}{ \times _{p \ne n}}\left\{ {{\mathbf \Omega} _p^{\left( n \right)}} \right\}$ where ${\mathbf \Omega}^{(n)}_p\in\mathbb{R}^{I_n\times R_n},\,\,p=1,2,\ldots,N$ are random matrices and then compute the unfolding matrix ${\mathbf W}_{(n)}$. In view of \eqref{TuckerUnfold}, the latter is equivalent to the multiplication of the unfolding matrix ${\mathbf W}_{(n)}$ with the Kronecker product of the random matrices\footnote{In \cite{sun2018tensor, che2018randomized}, the Khatri-Rao product of random matrices is used.}. This is a memory efficient approach because it needs $\mathcal{O}(\sum\limits_{p \ne n} {{I_p}{R_p}} )$ memory, while the former needs $\mathcal{O}(\left( {\prod\limits_{p \ne n} {{I_p}} } \right){R_n})$ memory. Also a memory efficient variant of the power iteration expression ${\left( {{{\mathbf Z}_{\left( n \right)}}{\mathbf Z}_{\left( n \right)}^T} \right)^q}{{\mathbf Z}_{\left( n \right)}}{{\mathbf \Omega}^{\left( n \right)}}$ in the tensor format is discussed in \cite{che2020computation}.

Both Algorithms \ref{RHOSVD} and \ref{HOOIRAND} are nonadaptive because they require an estimation of the multiliear rank $(R_1,R_2,\dots,R_N)$ which may be a difficult assumption in real applications.
To overcome this difficulty, it is possible to apply fixed precision randomized algorithms \cite{martinsson2016randomizedReveal, yu2018efficient} within Algorithms \ref{RHOSVD} and \ref{HOOIRAND} instead of standard (basic) randomized Algorithm \ref{ALG1}. This idea has been utilized in \cite{minster2019randomized}.

Also, a different strategy has been used in \cite{che2018randomized} for adaptively computing the multilinear rank of tensors. This idea basically relies on solving the Problems \eqref{DT1} and \eqref{DT2} numerically and finding orthonormal matrices ${\mathbf Q}^{(n)}$. 
Using the following identity
\[
\begin{array}{l}
\left\| \underline{\mathbf X} \right\|_F^2 = \left\| {\underline{\mathbf X}{ \times _n}\left( {{{\mathbf Q}^{(n)}}{\mathbf Q}^{(n)\,T}} \right)} \right\|_F^2+
\\
{\hspace{3.5cm}\left\| {\underline{\mathbf X}{ \times _n}\left( {{\mathbf I}_{I_n} - {{\mathbf Q}^{(n)}}{\mathbf Q}^{(n)\,T}} \right)} \right\|_F^2,}
\end{array}
\]
it can be shown that \cite{vannieuwenhoven2012new}
\begin{equation*}
\begin{array}{l}
\left\| {\underline{\mathbf X} - \underline{\mathbf X}{ \times _1}\left( {{{\mathbf Q}^{(1)}}{\mathbf Q}^{(1)\,T}} \right){ \times _2} \cdots { \times _N}\left( {{{\mathbf Q}^{(N)}}{\mathbf Q}^{N\,T}} \right)} \right\|_F^2 \\
\le \sum\limits_{n = 1}^N {\left\| {\underline{\mathbf X} - \underline{\mathbf X}{ \times _n}\left( {{{\mathbf Q}^{(n)}}{\mathbf Q}^{(n)\,T}} \right)} \right\|_F^2},
\end{array}
\end{equation*}
for $n=1,2,\ldots,N$.

This implies that finding orthonormal factor matrices is equivalent to finding the orthonormal matrices ${\mathbf Q}^{(n)}\in\mathbb{R}^{I_n\times R_n}$ that satisfy 
\begin{equation}\label{SubPro}
\begin{array}{l}
\left\| {\underline{\mathbf X} - \underline{\mathbf X}{ \times _n}\left( {{{\mathbf Q}^{(n)}}{\mathbf Q}^{(n)\,T}} \right)} \right\|_F^2 = \\ 
\hspace{2.75cm} \left\| {\underline{\mathbf X}{ \times _n}\left( {{\mathbf I}_{I_n} - {{\mathbf Q}^{(n)}}{\mathbf Q}^{(n)\,T}} \right)} \right\|_F^2 \le \frac{\varepsilon}{N}.
\end{array}
\end{equation}
The unfolding form of \eqref{SubPro} is 
\begin{equation}\label{SubProEQ}
\left\| {\left( {{{\mathbf I}_{{I_n}}} - {{\mathbf Q}^{(n)}}{\mathbf Q}^{(n)\,T}} \right){{\mathbf X}_{(n)}}} \right\|_F^2 \le \frac{\varepsilon}{N},
\end{equation}
and as a result, to find the factor matrix ${\mathbf Q}^{(n)}$, it is sufficient to compute an orthonormal matrix which captures the range of the unfolding matrix ${\mathbf X}_{(n)}$. Here, we can benefit from the efficient and optimized randomized algorithms for matrices. In \cite{che2018randomized}, the sub-problem \eqref{SubPro} which is in tensor form is solved.
\begin{center}
\RestyleAlgo{ruled}
\begin{algorithm}
\SetKwInOut{Input}{Input}
\SetKwInOut{Output}{Output}
	  \Input{A data tensor $\underline{\mathbf X} \in {\mathbb{R}^{{I_1} \times I_2\times \cdots  \times {I_N}}}$, and a multilinear rank $\left( {{R_1}, R_2, \ldots ,{R_N}} \right)$}
	        \Output{Approximative HOSVD of the tensor $\underline{\mathbf X}$ as 
	        $\underline{\mathbf X} \cong \left[\kern-0.15em\left[ {\underline{\mathbf S};{{\mathbf Q}^{(1)}},{{\mathbf Q}^{(2)}}, \ldots ,{{\mathbf Q}^{(N)}}}
 \right]\kern-0.15em\right]$}
      $\underline{\mathbf Z}=\underline{\mathbf X}$\\
%
  \For{$n = 1,2,\ldots, N$}
  {Compute ${{\mathbf W}^{\left( n \right)}} = {{\mathbf Z}_{(n)}}\,{{\boldsymbol\Omega} ^{\left( n \right)}},$ where ${{\boldsymbol\Omega}^{\left( n \right)}}$ is an $\left( {\prod\nolimits_{k \ne n} {{I_k}} } \right) \times {R_n}$ random Gaussian matrix\\
Compute ${\mathbf Q}^{(n)}$ as an orthonormal basis of ${\mathbf W}^{(n)}$ by using, e.g., the economic QR decomposition;\\
  }
    Compute the core tensor as $\underline{\mathbf S} \equiv \underline{\mathbf X}{ \times _1}{\mathbf Q}^{(1)\,T} { \times _2}{\mathbf Q}^{(2)\,T} \cdots { \times _N}{\mathbf Q}^{(N)\,T}$\\
 \caption{Random projection HOSVD (RP-HOSVD) Algorithm}\label{RHOSVD}
\end{algorithm}
\end{center}

\begin{center}
\RestyleAlgo{ruled}
\begin{algorithm}
\SetKwInOut{Input}{Input}
\SetKwInOut{Output}{Output}
	  \Input{A data tensor $\underline{\mathbf X} \in {\mathbb{R}^{{I_1} \times I_2\times \cdots  \times {I_N}}}$, and a multilinear rank $\left( {{R_1}, R_2, \ldots ,{R_N}} \right)$}
	        \Output{Approximative HOSVD of the tensor $\underline{\mathbf X}$ as
$\underline{\mathbf X} \cong \left[\kern-0.15em\left[ {\underline{\mathbf S};{{\mathbf Q}^{(1)}},{{\mathbf Q}^{(2)}}, \ldots ,{{\mathbf Q}^{(N)}}}
 \right]\kern-0.15em\right]$}
Initialize factor matrices ${\mathbf Q}^{(n)}\in\mathbb{R}^{I_n\times R_n}$ as factor matrices of the HOSVD or random Gaussian matrices;\hspace{5cm} \\
\While{A stopping criterion is not satisfied}{
  \For{$n = 1,2,\ldots, N$}
  {$\underline{\mathbf S} = \underline{\mathbf X}{ \times _{p \ne n}}\left\{ {{{\mathbf Q}^{\left( p \right)\,\,T}}} \right\}$\\
  Compute ${{\mathbf W}^{\left( n \right)}} = {{\mathbf S}_{(n)}}\,{{\boldsymbol \Omega} ^{\left( n \right)}}$ where ${{\boldsymbol \Omega} ^{\left( n \right)}} \in {\mathbb{R}^{\prod\limits_{p \ne n} {{R_p}}  \times {R_n}}}$ is a random matrix drawn from Gaussian distribution\\
Compute ${\mathbf Q}^{(n)}\in\mathbb{R}^{I_n\times R_n}$ as an orthonormal basis of ${\mathbf W}^{(n)}$, e.g., by using QR
decomposition
  }
  $\underline{\mathbf S}=\underline{\mathbf S}\times_{N} {\mathbf Q}^{(N)\,T}$
  }
 \caption{Random projection HOOI (RP-HOOI) algorithm}\label{HOOIRAND}
\end{algorithm}
\end{center}

Similar to the HOSVD and the HOOI algorithms, a randomized variant of the STHOSVD algorithm can also be developed. Algorithm \ref{RSeqHOSVD} is a random projection STHOSVD algorithm in which Algorithm \ref{ALG1} is exploited for computing low-rank approximations of the $n$-unfolding matrices \cite{minster2019randomized}. The adaptive randomized STHOSVD is also proposed in \cite{minster2019randomized} wherein the randomized rank-revealing algorithms  \cite{yu2018efficient,martinsson2016randomizedReveal} are exploited instead of the Algorithm \ref{ALG1}. All discussed algorithms so far do not preserve the structure of the underlying data tensor such as sparsity or nonnegativity. A structural preserving variant of Algorithm \ref{RSeqHOSVD} is developed in \cite{minster2019randomized} where the core tensor takes elements from the original data matrix.

\begin{center}
\RestyleAlgo{ruled}
\begin{algorithm}
\SetKwInOut{Input}{Input}
\SetKwInOut{Output}{Output}
	  \Input{A data tensor $\underline{\mathbf X} \in {\mathbb{R}^{{I_1} \times I_2\times \cdots  \times {I_N}}}$, and a multilinear rank $\left( {{R_1}, R_2, \ldots ,{R_N}} \right)$}
	        \Output{Approximative Tucker representation of the tensor $\underline{\mathbf X}$ as
$\underline{\mathbf X} \cong \left[\kern-0.15em\left[ {\underline{\mathbf S};{{\mathbf Q}^{(1)}},{{\mathbf Q}^{(2)}}, \ldots ,{{\mathbf Q}^{(N)}}}
 \right]\kern-0.15em\right]$}
${\mathbf S}={\mathbf X}_{(n)}$\\
  \For{$n = 1,2,\ldots, N$}
  {
  $[{\mathbf Q}^{(n)},\mathtt{\sim},\mathtt{\sim}]=$  Apply Algorithm \ref{ALG1} to the $n$-unfolding matrix ${\mathbf S}_{\left(n\right)}$ with target rank $R_n$\\
$\underline{\mathbf S}=\underline{\mathbf S}\times_{n} {\mathbf Q}^{(n)\,T}$\\
  }
 \caption{Randomized Sequentially Truncated HOSVD (R-STHOSVD) algorithm}\label{RSeqHOSVD}
\end{algorithm}
\end{center}

The algorithms and techniques discussed so far need to pass the original data tensor $\underline{\mathbf X}$ multiple times. This restricts their applicability for the data tensors stored out-of-cores because of high communication cost. Due to this issue, we need to develop so-called randomized Pass-Efficient algorithms where the original data tensor should be passed as few as possible. For example, one possible option is using Algorithm \eqref{TwoSide-SVD} or pass-efficient algorithms proposed in \cite{halko2011finding, woodruff2014sketching, boutsidis2016optimal, upadhyay2018price, bajaj2019sketchycoresvd}, for low-rank approximation of $n$-unfolding matrices.

An elegant pass-efficient algorithm is developed in \cite{sun2019low} and is described in Algorithm \ref{Onestreaming}. It first captures the important actions of the tensor $\underline{\mathbf X}$ in each mode (Lines 1-2), i.e., the range of unfolding matrices, and also their corresponding interactions (Line 3), i.e., the compressed core tensor generated by multiplying the original data tensor with random matrices in different modes. These are performed by multiplying the unfolding matrices by random matrices and at the same time multiplying an original tensor by random matrices along its different modes. The structure of this procedure is described in Lines 1 and 3 of Algorithm \ref{Onestreaming} where
\[
{{\mathbf \Omega}_n} \in {\mathbb{R}^{\left(\prod\limits_{i \ne n} {{I_i}}\right) \times {K_n}}},\,\,\,\,\widetilde{{\mathbf \Omega}}_n \in {\mathbb{R}^{{S_n} \times {I_n}}},
\]
and $R_n \le K_n \le S_n,\,\,n=1,2,\ldots,N.$
Here, $(R_1,R_2,\ldots,R_N)$ is a target multilinear rank.
Note that in Algorithm \ref{Onestreaming}, ${{\mathbf Q}^{\left( n \right)}} \in {\mathbb{R}^{{I_n} \times {K_n}}},\,\,\,\,\underline{\mathbf H} \in {\mathbb{R}^{{S_1} \times {S_2}\times  \cdots  \times {S_N}}}$, and in order to obtain a truncated HOSVD, the obtained core tensor $\underline{\mathbf S}$ and its corresponding factor matrices ${\mathbf Q}^{(N)}$ are truncated to multilinear rank $(R_1,R_2,\ldots,R_N)$.

\begin{center}
\RestyleAlgo{ruled}
\begin{algorithm}
\LinesNumbered
\SetKwInOut{Input}{Input}
\SetKwInOut{Output}{Output}
	  \Input{A data tensor $\underline{\mathbf X}\in\mathbb{R}^{I_1\times I_2\times \cdots \times I_N}$, sketching parameters $K,\,S$, $S > K$ and a predefined multilinear rank $(R_1,R_2,\ldots,R_N)$}
	        \Output{Approximative Tucker representation of the tensor $\underline{\mathbf X}$ as $\underline{\mathbf X} \cong \left[\kern-0.15em\left[ {\underline{\mathbf S};{{\mathbf Q}^{(1)}},{{\mathbf Q}^{(2)}}, \ldots ,{{\mathbf Q}^{(N)}}}
 \right]\kern-0.15em\right]$}
Compute factor sketches, ${{\mathbf Y}_n} = {{\mathbf X}_{\left( n \right)}}\,{{\mathbf \Omega}_n},\,n = 1, 2,\ldots $ $,N$\\
Recover factor matrices, $\left[ {{{\mathbf Q}^{(n)}},\mathtt{\sim}} \right] = {\mathbf Q}{\mathbf R}\left( {{{\mathbf Y}_n}} \right)$\\
Compute core sketch $\underline{\mathbf H} = \underline{\mathbf X}{ \times _1}\widetilde{\mathbf \Omega}_1\, \times\, \cdots { \times _N}\,\widetilde{\mathbf \Omega}_N$\\
Recover core tensor $\underline{\mathbf S} = \underline{\mathbf H}{ \times _1}{\left( {{\widetilde{\mathbf \Omega}}_1{{\mathbf Q}^{\left( 1 \right)}}} \right)^\dag }\times_2 \cdots { \times _N}\,{\left( {{\widetilde{\mathbf \Omega}}_N{{\mathbf Q}^{\left( N \right)}}} \right)^\dag }$\\
Truncate the core tensor $\underline{\mathbf S}$ and the factor matrices ${\mathbf Q}^{(n)},\,n=1,2,\ldots,N$ with multilinear rank $(R_1,R_2,\ldots,R_N)$
 \caption{Randomized pass-efficient randomized algorithm for the Tucker decomposition (R-PET)}\label{Onestreaming}
\end{algorithm}
\end{center}

\subsection{Randomized Sampling Tucker Decomposition}\label{RSamp}
In randomized sampling Tucker decomposition, each factor matrix ${\mathbf Q}^{(n)}$ is computed by sampling the columns of the corresponding unfolding matrix ${\mathbf X}_{(n)}$ or equivalently sampling the fibers\footnote{Algorithms based on sampling slices are proposed in \cite{tarzanagh2018fast, mahoney2008tensor, song2019relative}, but those in \cite{tarzanagh2018fast, mahoney2008tensor} are not related to the Tucker decomposition.} of the original data tensor $\underline{\mathbf X}$. The structure of this approach is outlined in Algorithm \ref{SamTucker} \cite{drineas2007randomized}. After computing an approximation for the range of unfolding matrices, they can be used to find a Tucker approximation \eqref{TK1}, although they can also be orthogonalized by economic QR decomposition to be used for computation of the HOSVD. Note that in the R-ST algorithm, the factor matrices preserve the structure of the original data tensor (nonnegativity, sparsity, smoothness) but not  necessarily the core tensor. As we mentioned in Section \ref{Sec:RS}, the procedure of column or row selection can be performed based on different kinds of probability distributions such as uniform or other distributions and also with or without replacement. Moreover, the algorithms naturally provide either additive or relative errors.
Sampling based on the leverage scores are proposed in \cite{cheng2016spals, perros2015sparse, saibaba2016hoid}.

\begin{center}
\RestyleAlgo{ruled}
\begin{algorithm}
\SetKwInOut{Input}{Input}
\SetKwInOut{Output}{Output}
	  \Input{A data tensor $\underline{\mathbf X}\in \mathbb{R}^{I_1\times I_2\times \cdots\times I_N}$,\,and a predefined multilinear rank $(R_1,R_2,\ldots,R_N)$}
	        \Output{Approximative Tucker representation of the tensor $\underline{\mathbf X}$ as
$\underline{\mathbf X} \cong \left[\kern-0.15em\left[ {\underline{\mathbf S};{{\mathbf Q}^{(1)}},{{\mathbf Q}^{(2)}}, \ldots ,{{\mathbf Q}^{(N)}}}
 \right]\kern-0.15em\right]$}
\For{$n=1,2,\ldots,N$}
{
Sample some columns of ${\mathbf X}_{(n)}$ based on a probability distribution and store them in the factor matrix ${\mathbf Q}^{(n)}\in \mathbb{R}^{I_n\times R_n}$
}
Compute the core tensor $\underline{\mathbf S}=\underline{\mathbf X}{ \times _1}{{\mathbf Q}^{\dag}_1}{ \times _2}{{\mathbf Q}^{\dag}_2} \cdots { \times _N}{{\mathbf Q}^{\dag}_N}$ 
 \caption{Randomized Sampling Tucker approximation (R-ST)}\label{SamTucker}
\end{algorithm}
\end{center}

Cross-approximation\footnote{It is also called (pseudo)-skeleton decomposition or CUR decomposition.} 
\cite{gantmacher1959theory, goreinov1997theory, goreinov1997pseudo}, can be considered as a sampling technique with a difference that the sampling procedure is performed heuristically instead of randomly. 
Three known heuristic algorithms are: 
\begin{itemize}
\item Maxvol-based low-rank matrix approximation \cite{goreinov2010find,goreinov2001maximal}

\item Cross2D matrix approximation \cite{savostyanov2006polilinear,tyrtyshnikov2000incomplete}

\item Discrete Empirical Interpolatory Method (DEIM) \cite{chaturantabut2009discrete,sorensen2016deim}

\item Pivoted QR decomposition \cite{stewart1999four}
\end{itemize}
It is know that the quality of the cross approximation quite depends on the module of the determinant of the intersection matrix\footnote{By an intersection matrix, we mean a matrix which produced by the intersection of sampled columns and rows.} which is called {\em matrix volume}. 
More precisely, a set of columns and rows should be selected with an intersection matrix whose volume is as much as possible. Clearly, this is an NP hard problem because we need to check the volume of all possible intersection matrices produced by different selection of columns and rows. However, heuristic algorithms exist for computing suboptimal solutions. 

Several generalizations of the matrix cross-approximation to the low Tucker rank approximation are proposed in \cite{caiafa2010generalizing,oseledets2008tucker}. The proposed algorithms in \cite{oseledets2008tucker} basically apply the Cross2D algorithm \cite{savostyanov2006polilinear,tyrtyshnikov2000incomplete} to the unfolding matrices sequentially where the long rows of the unfolding matrices can be treated as slice matrices. Here again a cross approximation of these slices are computed and this idea totally reduces the computational complexity of the cross-approximation. Two algorithms proposed in \cite{caiafa2010generalizing} are based on the fiber selection and they do not apply cross-approximation to unfolding matrices. The first algorithm is analytic while second one is iterative (adaptive). The analytic one is of less practical interest because it requires quite large number of fibers, but the adaptive one is more efficient and is a straightforward generalization of cross-approximation of matrices to tensors. A similar
approach is proposed in \cite{friedland2011fast} in the sense of number of selected
fibers in each mode.

The fibers can also be sampled using the pivoted QR decomposition applied to the unfolding matrices ${\mathbf X}_{(n)}$. Inspired by the matrix case \cite{stewart1999four}, {\em Higher Order Interpolatory Decomposition (HOID)} is proposed in \cite{saibaba2016hoid}. Given an $N$th-order tensor $\underline{\mathbf X}\in\mathbb{R}^{I_1\times I_2\times\cdots\times I_N}$, the pivoted QR factorization is applied to the unfolding matrix ${\mathbf X}_{(n)} \in {\mathbb{R}^{{I}_1 \times {J}}},\,(J = \prod\limits_{i \ne n} {{I_i}})$ as
\begin{equation}\label{PQR}
{\mathbf X}_{(n)}{\mathbf \Pi}  = {\mathbf Q} {\mathbf R},
\end{equation}
where ${\mathbf\Pi}\in{\mathbb{R}^{J \times J}}$ and ${\mathbf Q} \in {\mathbb{R}^{{I}_n \times {I}_n}}$ are a permutation matrix and an orthogonal matrix, respectively. The pivoted QR decomposition can be computed using strong rank-revealing QR (RRQR) algorithm \cite{Gu}. The permutation matrix ${\mathbf \Pi}$ and the corresponding orthogonal matrix $\mathbf Q$ are partitioned as follows
\begin{equation}\label{Permu}
{\mathbf \Pi}=\left[ {{{\mathbf \Pi} _1}\,\,{{\mathbf \Pi} _2}} \right],\quad\quad  {\mathbf Q}=\left[ {{\mathbf Q}_1\,\,{\mathbf Q}_2} \right],
\end{equation}
where ${\mathbf\Pi}_1\in {\mathbb{R}^{{J} \times K}},\,\,\,\,{{\mathbf\Pi}_2} \in {\mathbb{R}^{{J} \times \left( {J - K} \right)}},\,\,\,{{\mathbf Q}_1} \in {\mathbb{R}^{{I}_n \times K}},\,\,\,\,\,{{\mathbf Q}_2} \in {\mathbb{R}^{{I}_n \times \left( {I_n - K} \right)}}$.
Here Equation \eqref{PQR} can be rewritten as
\begin{equation}\label{PQR2}
{\mathbf X}_{(n)}\left[ {{{\mathbf\Pi}_1}\,\,{{\mathbf\Pi}_2}} \right] = \left[ {{{\mathbf Q}_1}\,\,{{\mathbf Q}_2}} \right]\left[ {\begin{array}{*{20}{c}}
{{{\mathbf R}_{11}}}&{{{\mathbf R}_{12}}}\\
{\mathbf 0}&{{{\mathbf R}_{22}}}
\end{array}} \right],
\end{equation}
 where ${\mathbf R}_{11} \in \mathbb{R}^{K \times K},\,\,\,{\mathbf R}_{12} \in \mathbb{R}^{K \times \left( {I_n-K} \right)},\,\,{\mathbf R}_{22} \in {\mathbb{R}^{\left( {I_n - K} \right) \times \left( {J - K} \right)}}$ and ${\mathbf 0}\in\mathbb{R}^{(I_n-K)\times K}$ with all entries equal to zero.
From \eqref{PQR2}, by straightforward computations, we have 
\[
\begin{array}{l}
{\mathbf Q}^{(n)} \equiv {\mathbf X}_{(n)}{{\mathbf \Pi} _1} = {{\mathbf Q}_1}{{\mathbf R}_{11}},\\
\\
{\mathbf X}_{(n)}{{\mathbf\Pi}_2} = {{\mathbf Q}_1}{{\mathbf R}_{12}} + {{\mathbf Q}_2}{{\mathbf R}_{22}} \cong {{\mathbf Q}_1}{{\mathbf R}_{12}},
\end{array}
\]
if $\left\| {\mathbf R}_{22} \right\|_2$ is small enough. It can be seen that 
\begin{equation}\label{Ap}
{\mathbf X}_{(n)} \cong \left[ {{{\mathbf Q}_1}{{\mathbf R}_{11}},\,\,{{\mathbf Q}_1}{{\mathbf R}_{12}}} \right]{{\mathbf\Pi}^T},
\end{equation}
and substituting ${{\mathbf Q}_1} = {\mathbf Q}^{(n)}{\mathbf R}_{11}^{-1}$ in \eqref{Ap}, we have
\begin{equation}\label{InterpoFormul}
{\mathbf X}_{(n)} \cong {\mathbf Q}^{(n)}{{\mathbf F}^T},\,\,\,\,\,{\mathbf F}^{T} =\left[ {{\mathbf I}\,\,\,\,{\mathbf R}_{11}^{ - 1}{{\mathbf R}_{12}}} \right]{{\mathbf \Pi}^T},
\end{equation}
which is a low-rank matrix approximation of the matrix ${\mathbf X}_{(n)}$. The matrix ${\mathbf Q}^{(n)}$ which is a full-rank matrix\footnote{Because it is the multiplication of a nonsingular matrix and an orthogonal matrix.} can be used as an approximation for the basis of the range of ${\mathbf X}_{(n)}$. It is worth mentioning that ${\mathbf Q}^{(n)}={\mathbf X}_{(n)}(:,{\mathbf p})$ where ${\mathbf p}$ is a set indices of the selected columns and as a result it may not be necessarily orthonormal. The HOID algorithm applies the earlier procedure to all unfolding matrices ${\mathbf X}_{(n)}$ and computes the basis matrices ${\mathbf Q}^{(n)}$. Afterwards, the core tensor $\underline{\mathbf S}$ can be computed as follows 
 \[
  \underline{\mathbf S} = \underline{\mathbf X}{ \times _1}{\mathbf Q}^{(1)\,\dag}\times_2 {\mathbf Q}^{(2)\,\dag} \cdots { \times _N}{\mathbf Q}^{N\, \dag}.
  \]
A randomized variant of this algorithm can be developed by replacing randomized QR decomposition instead of the QR decomposition. To this end, we should first make reduction in the second mode of the unfolding matrices using random projection as ${\mathbf Y} = {\mathbf X}_{(n)}{\mathbf \Omega}$ where ${\mathbf\Omega}  \in {\mathbb{R}^{\prod\limits_{i \ne n} {{I_n}}\times R}}$ is a random matrix after which the economic QR decomposition of the matrix ${\mathbf Y}$ is computed as ${\mathbf Y}={\mathbf Q}{\mathbf R}$. The orthonormal matrix ${\mathbf Q}$ is an approximate orthonormal basis for the range of the matrix ${\mathbf X}_{(n)}$ which can be used in the above procedure. This procedure is summarized in Algorithm \ref{HOID}. In situation that the singular values of the matrix ${\mathbf X}_{(n)}$ do not decay very fast, the orthonormal basis may not be accurate and power iteration technique is required (see Remark \ref{rempo}).

If the dimension of the first mode is also large, it is suggested in \cite{voronin2017efficient} to first make reduction in the first mode of the matrix using random projection as ${\mathbf Y} = {\mathbf \Omega}{\mathbf X}_{(n)}$ where ${\mathbf\Omega}  \in {\mathbb{R}^{(R+P)\times I_n}}$ is a random matrix ($P$ is the oversampling parameter), after which the procedure described above is applied to the matrix ${\mathbf Y}$ for column selection. It is shown that the selected indices of columns for the matrix ${\mathbf Y}$ can be used for the original data matrix ${\mathbf X}$. More precisely, if ${\mathbf Y} \cong {\mathbf Y}({:,{\mathbf p}}){\mathbf F}^T$ then the indices ${\mathbf p}$ corresponding to the selected columns can be used for the matrix ${\mathbf X}_{(n)}$, i.e., ${\mathbf X}_{(n)}\cong {\mathbf X}_{(n)}(:,{\mathbf p}){\mathbf F}^T$.


\begin{rem}
The procedure of column selection can be performed using the leverage scores sampling or the DEIM algorithm. These techniques were studied in \cite{saibaba2016hoid}.
\end{rem}

\begin{center}
\RestyleAlgo{ruled}
\begin{algorithm}
\SetKwInOut{Input}{Input}
\SetKwInOut{Output}{Output}
	  \Input{A data tensor $\underline{\mathbf X} \in {\mathbb{R}^{{I_1} \times I_2\times \cdots  \times {I_N}}}$, and a multilinear rank $\left( {{R_1}, R_2, \ldots ,{R_N}} \right)$}
	        \Output{Approximative Tucker representation of the tensor $\underline{\mathbf X}$ as $\left[\kern-0.15em\left[ {{\underline{\mathbf S},\,\,{{\mathbf Q}^{\left( 1 \right)}}, {{\mathbf Q}^{\left( 2 \right)}}, \ldots ,{{\mathbf Q}^{\left( N \right)}}}} 
 \right]\kern-0.15em\right]$}
  \For{$n = 1,2,\ldots, N$}
  {Compute a randomized interpolatory decomposition of unfolding ${{\mathbf X}_{\left( n \right)}} \cong {{\mathbf Q}^{(n)}}{\mathbf F}_n^T,$
 where $\,\,\,{{\mathbf Q}^{(n)}} \in {\mathbb{R}^{{I_n} \times {R_n}}}$, are columns of $\,\,{{\mathbf X}_{(n)}} 
$
  }
  Compute the core tensor $\underline{\mathbf S} \in {\mathbb{R}^{{R_1} \times {R_2}\times \cdots  \times {R_N}}}$ as 
  \[
  \underline{\mathbf S} \equiv \underline{\mathbf X}{ \times _1}{\mathbf Q}^{(1)\,\dag}\times_2  \cdots { \times _N}{\mathbf Q}^{N\, \dag}
  \]
 \caption{Randomized Higher Order Interpolatory Decomposition (HOID) (R-HOID) Algorithm}\label{HOID}
\end{algorithm}
\end{center}
Other sampling techniques for the HOSVD and the HOOI algorithms can be found in \cite{tsourakakis2010mach, song2019relative, traore2019singleshot}.
The proposed algorithms in \cite{tsourakakis2010mach} are based on the sparsification idea where in the first step, a sparse tensor is generated from the original data tensor. Then, it is used for subsequent computations. The idea of the stochastic gradient descent is used in \cite{traore2019singleshot}, where at each iteration only a subtensor of the original tensor is considered. The randomized algorithms proposed in \cite{song2019relative} are based on sampling fibers instead of components.

\subsection{Randomized least-squares Tucker Decomposition}
At each iteration of the randomized Algorithm \ref{HOOIRAND}, several multiplications with random matrices and also contraction of the core tensor with factor matrices are required which may be expensive if the algorithm needs many iterations to converge. Motivated by this fact, in \cite{malik2018low}, a randomized least-squares HOOI algorithm is proposed. This algorithm is summarized in Algorithm \ref{TUCKER-TS} in which instead of exploiting updating rules 4 and 7 in Algorithm \ref{AHOOI}, the equivalent least-squares problems \eqref{LSS} and \eqref{LSE2} are solved, respectively. To be precise, instead of performing Lines 4-5 and 7 of Algorithm \ref{AHOOI}, the least-squares problem \eqref{LSS} is solved by random least-squares algorithms. Please note that the coefficient matrices in least-squares problems \eqref{LSS} and \eqref{LSE2} are not required to be computed explicitly. Also, the TensorSketch \cite{pagh2013compressed} is used in \cite{malik2018low} to solve the underlying least-squares problem. 

\begin{center}
\RestyleAlgo{ruled}
\begin{algorithm}
\LinesNumbered
\SetKwInOut{Input}{Input}
\SetKwInOut{Output}{Output}
	  \Input{A data tensor $\underline{\mathbf X} \in {\mathbb{R}^{{I_1} \times I_2\times \cdots  \times {I_N}}}$, and a multilinear rank $\left( {{R_1}, R_2, \ldots ,{R_N}} \right)$;}
	        \Output{Approximative Tucker representation of the tensor $\underline{\mathbf X}$ as $\underline{\mathbf X} \cong \left[\kern-0.15em\left[ {\underline{\mathbf S};{{\mathbf Q}^{(1)}},{{\mathbf Q}^{(2)}}, \ldots ,{{\mathbf Q}^{(N)}}}
 \right]\kern-0.15em\right]$}
	        Initialize factor matrices ${\mathbf Q}^{(1)},{\mathbf Q}^{(2)},\ldots,{\mathbf Q}^{(N)}$ and core tensor $\underline{\mathbf S}$ with i.i.d uniformly random entries\\
	        \While{A stopping criterion is not satisfied}{
  \For{$n = 1,2,\ldots, N$}
  {  Compute ${\mathbf Q}^{(n)}$ by solving least-squares problem \eqref{LSS} via a sketching technique for $\,\,n=1,2,\ldots,N$;
  }
  {Find the core tensor $\underline{\mathbf S}$ by solving least-squares problem \eqref{LSE2} via sketching technique;
  }
  }
 \caption{Randomized least-squares Higher Order Orthogonal Iteration (HOOI) (R-LSHOOI) Algorithm}\label{TUCKER-TS}
\end{algorithm}
\end{center}

\subsection{Randomized Count-Sketch Tucker Decomposition}
The count-sketch technique discussed in Subsection \ref{count-sketch}, can be straightforwardly generalized to tensors. Motivated by the paper \cite{battaglino2018practical} in which the underlying least-squares problems in the CP-ALS\footnote{ALS means Alternating Least-Squares \cite{kolda2009tensor}.} algorithm \cite{kolda2009tensor} are solved via the sampling techniques, in \cite{malik2018low}, the count-sketch technique is used to solve the least-squares problems \eqref{LSS} and \eqref{LSE2}. More precisely, the reduction map ${\mathbf T}$ in the least-squares problem \eqref{RLST}, is the count-sketch operator. Because of a special structure of the coefficient matrix of least-squares problems \eqref{LSS} and \eqref{LSE2}, sketching procedure can be performed very fast by employing FFT \cite{malik2018low}. 

The count-sketch technique is used in \cite{malik2019fast} for tensor interpolative decomposition \cite{biagioni2015randomized} and also in \cite{wang2015fast} for the CP decomposition. The concept of Higher-order Count Sketch is developed in \cite{shi2019multi} for higher order tensors to fully exploit the multidimensional structure of the data tensors. 

Another possible direction is using the count-sketch technique for low-rank matrix approximations of the unfolding matrices which are required for tensor decomposition. Here, the factor matrices are approximated by computing low-rank approximation of the unfolding matrices using the count-sketch technique and then the core tensor is computed via \eqref{Core2}.

\section{\bf Application of the randomized HOSVD in fast Canonical Polyadic Decomposition (CPD)}\label{Sec:Appl}
The idea of making a prior reduction in a raw data tensor into the Tucker format and then decomposing its core tensor in the CPD format, was first proposed by Bro in his PhD thesis \cite{bro1998multi}. 
This technique is applicable only if the rank of the tensor does not exceed its dimensions. 
Motivated by this difficulty, a prior reduction in the Tensor Train (TT) format \cite {oseledets2011tensor} was recently suggested in \cite{phan2018tensor}. In this section, we focus on the former and discuss how randomization can be exploited within 
the tensor decomposition procedure. 

The main procedure of this technique is as follows
\begin{itemize}
\item Dimensionally reduction of a given data tensor using the Tucker decomposition as a preprocessing step,

\item Computing the CPD of the compressed core tensor of the Tucker decomposition,

\item Recovering the factor matrices of the CPD of the original data tensor from the factor matrices of the CPD of the core tensor.
\end{itemize}

Assuming that the data tensor has a low multilinear rank approximation, the first step mentioned above can be done through the randomized HOSVD algorithms.
This makes the decomposition procedure much faster compared with the deterministic ones. This idea is used in \cite{erichson2017randomized}.

Assume that the Tucker representation of a data tensor $\underline{\mathbf X}$ is computed as $\left[\kern-0.15em\left[ {\underline{\mathbf S};{{\mathbf Q}^{\left( 1 \right)}},{{\mathbf Q}^{\left( 2 \right)}}, \ldots ,{{\mathbf Q}^{\left( N \right)}}} \right]\kern-0.15em\right]$ using one of the randomized HOSVD or randomized Tucker decomposition algorithms discussed in Section \ref{Sec:RHOSVD}. Let $\underline{\mathbf S}\cong\left[\kern-0.15em\left[ {{\widetilde{\mathbf A}^{\left( 1 \right)}},{\widetilde{\mathbf A}^{\left( 2 \right)}}, \ldots ,{\widetilde{\mathbf A}^{\left( N \right)}}}\right]\kern-0.15em\right]$ be the CPD of the compressed core tensor $\underline{\mathbf S}$. Then, the CPD of the original tensor $\underline{\mathbf X}$ can be recovered from the CPD of the compressed tensor $\underline{\mathbf S}$. More precisely, the CPD of the original tensor $\underline{\mathbf X}$ is 
\[\underline{\mathbf X}\cong\left[\kern-0.15em\left[ {{{\mathbf A}^{\left( 1 \right)}},{{\mathbf A}^{\left( 2 \right)}}, \ldots ,{{\mathbf A}^{\left( N \right)}}} 
 \right]\kern-0.15em\right],\]
where 
${{{\mathbf A}^{\left(n\right)}}}={\mathbf Q}^{(n)}{{\widetilde{\mathbf A}^{\left(n\right)}}},\,n=1,2,\ldots,N.$
\begin{rem}
The idea of a prior reduction of a data tensor into the TT format for computation of the CPD was suggested in \cite{phan2018tensor}. Randomized variants of these algorithms can also be developed.
\end{rem}

\section{Discussion on further challenges}\label{Sec:FuthChall}
Recently, several scalable algorithms were proposed in \cite{oh2018scalable, oh2019high, oh2017s} for computation of the HOSVD.
These algorithms inherently do not have a randomized structure and basically they exploit the idea of {\em on-the-fly-computation} and {\em parallel row-wise update rule roles} to avoid the {\em intermediate data explosion} problem. Combining these algorithms with randomized techniques is a potential topic needs to be investigated. For example, this issue was recently studied in \cite{jang2020d}. Block Term Decomposition (BTD) \cite{de2008decompositionsI, de2008decompositionsII, de2008decompositionsIII} is a generalization of the CP decomposition and the Tucker decomposition and has found many applications 
such as blind source separation \cite{de2012block}, feature extraction \cite{ribeiro2016enhanced}, electroencephalogram (EEG) analysis \cite{hunyadi2014block} etc. Proposing randomized algorithms for this tensor decomposition is a potential research topic needs to be investigated.

Motivated by some applications in Cyber-Physical-Social Systems and Internet of Things (IOT), distributed algorithms for computation of the HOSVD were recently developed in \cite{yang2018multi, wang2018distributed, wang2018improved, wang2017big, yang2017tensor}. The proposed algorithms mainly distribute the unfolding matrices among several machines and integrate their low-rank matrix approximations to find the HOSVD approximation of the original data tensor. These algorithm can be further improved by exploiting randomized algorithms. 

\section{Simulations}\label{Sec:Sim}
In this section, we experimentally examine some of the selected and most efficient randomized algorithms presented in Section \ref{Sec:RHOSVD} and compare their performance and efficiency. All numerical simulations were performed on a cluster with 508Gb RAM and 48 CPUs 1.2GHz.
Although, the sampling procedure can be done based on different types of probability distributions, but we used uniform distribution in our computations. Note that in the case of randomized sampling algorithms, we have used sampling without replacement as it works better than sampling with replacement in machine learning applications \cite{kumar2009sampling, drinea2001randomized}. We experienced the same results. We have set oversampling $P=10$ and power iteration $q=2$ for all random projection algorithms (expect the noiseless case where we have not used any oversampling or power iteration parameter). For the R-PET algorithm (Algorithm \ref{Onestreaming}), we have used $S=2K+1$ in all experiments where for each example, a specific $K$ was used and for the STHOSVD algorithm, we always considered the ascending order $p=[1,2,\ldots,N]$. The relative error is used as a criterion to assess the performance of the randomized algorithms for the Tucker approximation and is defined as follows 
\begin{equation}\label{Relativeerror}
E = \frac{{{{\left\| {\underline{\mathbf X} - \widehat{\underline{\mathbf X}}} \right\|}_F}}}{{{{\left\| \underline{\mathbf X} \right\|}_F}}},
\end{equation}
where $\underline{\mathbf X}$ is a given data tensor and $\widehat{\underline{\mathbf X}}$ is its Tucker approximation, respectively.  

We also consider the case that the tensor $\underline{\mathbf X}$ is corrupted by a noise term as follows
\begin{equation}\label{noisytensor}
\underline{\mathbf Y} = \underline{\mathbf X} + \gamma\underline{\mathbf N},
\end{equation}
where $\underline{\mathbf N}$ is a Gaussian random tensor\footnote{A tensor whose elements are independent and identically distributed (i.i.d) and taken from normal distribution with zero mean and variance 1.} and $\gamma$ is a parameter to change the noise level. The Signal-to-Noise Ratio (SNR) measure is used to evaluate the noise level and is defined as follows 
\[
{\rm{SNR}}\,\left[ {{\rm{dB}}} \right] = 10\,\log\left( {\frac{{\left\| \underline{\mathbf X} \right\|_F^2}}{{\left\| {\gamma \underline{\mathbf N}} \right\|_F^2}}} \right).
\]
For the noisy case, the performance of an algorithm is evaluated in terms of fit which is defined as follows
\[
{\rm Fit} = 1 - E,
\]
where $E$ was defined by \eqref{Relativeerror}. The compression ratio is defined as the ratio $\frac{S_1}{S_2}$ where $S_1$ is the number of components of a given data tensor ${\underline{\mathbf X}}$ and $S_2$ is the number of components of the factor matrices and the core tensor obtained from the Tucker decomposition of ${\underline{\mathbf X}}$.


\subsection{Test 1 -- synthetic data tensors} 
The synthetic data tensors used in our experiments are: {random tensors with low multilinear rank} (noiseless and noisy cases), {random sparse tensors}, {Hilbert tensors} and {function based tensors}. 
\subsubsection{\bf Random Low Multilinear Rank Tensors}
\begin{figure}[!h]
\includegraphics[width=1\linewidth]{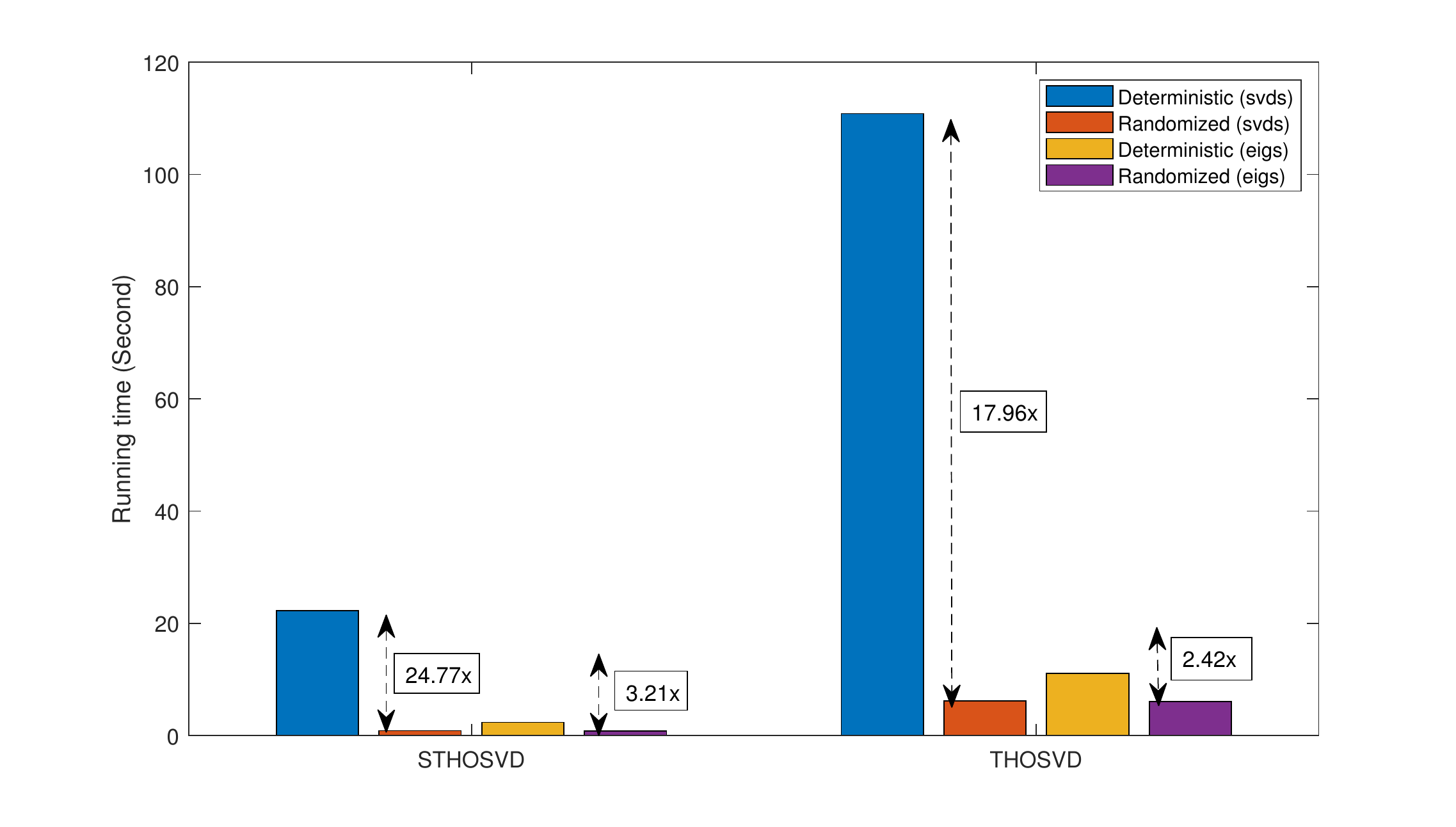}
\centering
\caption{\small Running time comparison of deterministic and randomized variants of the THOSVD and the STHOSVD algorithms (svds and eigs versions) for a noiseless random data tensor of size $1000\times 1000\times 1000\,(10^9\,\,{\rm enteries})$ and multilinear rank $(20,40,30)$.}
\label{FIGrand1}
\end{figure}

\begin{figure}[!h]
\includegraphics[width=1\linewidth]{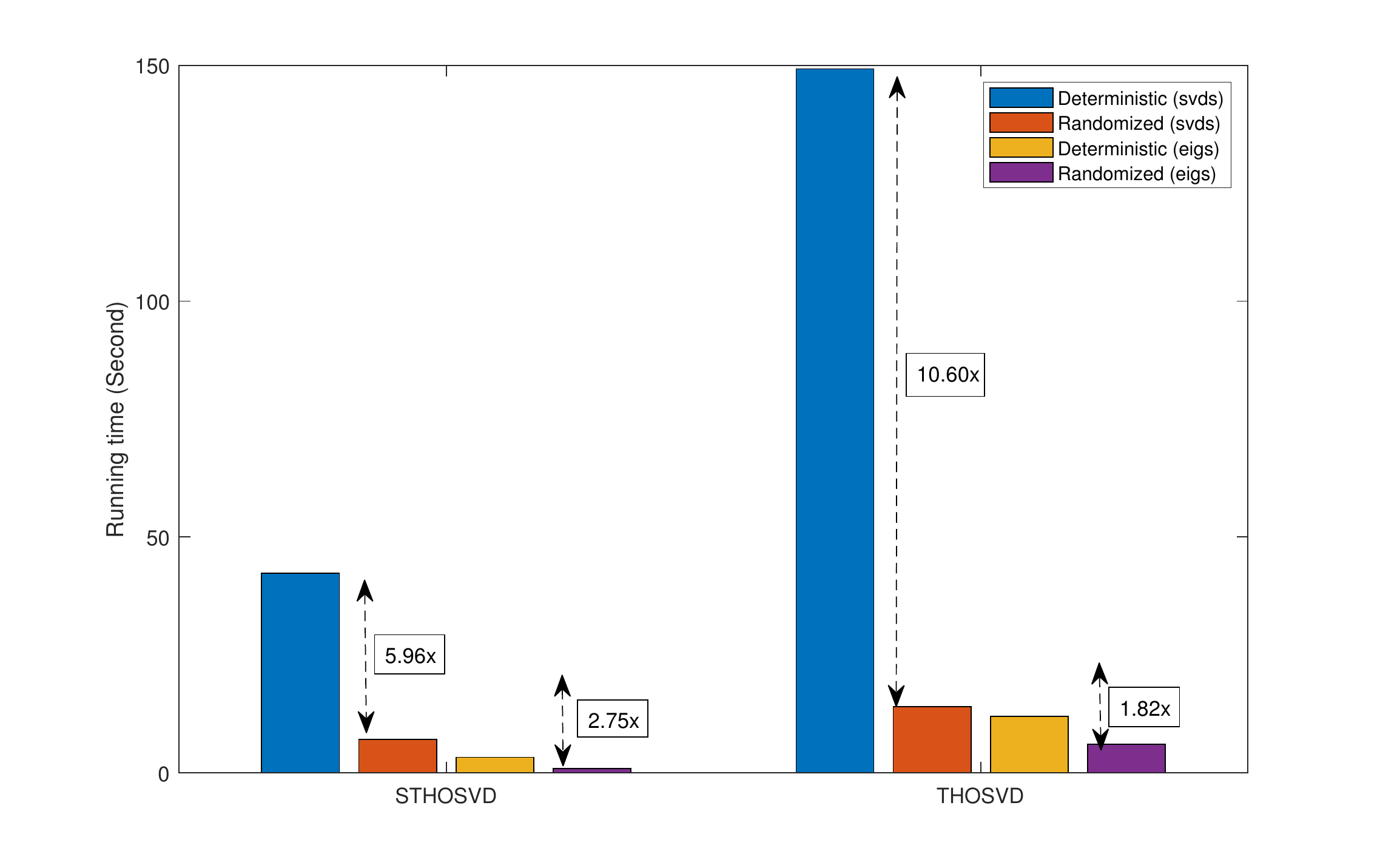}
\centering
\caption{\small Running time comparison of deterministic and randomized variants of the THOSVD and the STHOSVD algorithms (svds and eigs versions) for a noisy random data tensor of size $1000\times 1000\times 1000$ and multilinear rank $(20,40,30)$ with SNR=20 dB.}
\label{FIGrand2}
\end{figure}

In this experiment, we consider a 3rd order tensor of size $1000\times 1000\times 1000$ with multilinear rank $(20,40,30)$. We first randomly generate a Gaussian core tensor $\underline{\mathbf S}\in\mathbb{R}^{20\times 40\times 30}$ and Gaussian factor matrices ${\mathbf Q}^{(1)}\in\mathbb{R}^{1000\times 20},\,{\mathbf Q}^{(2)}\in\mathbb{R}^{1000\times 40}\,\,{\mathbf Q}^{(3)}\in\mathbb{R}^{1000\times 30}$. Then, the following tensor 
\begin{equation}\label{noiselesstensor}
\underline{\mathbf X} = \left[\kern-0.15em\left[ {\underline{\mathbf S};{{\mathbf Q}^{(1)}},{{\mathbf Q}^{(2)}},{{\mathbf Q}^{(3)}}}
 \right]\kern-0.15em\right],
\end{equation}  
is generated and some of the deterministic and randomized algorithms are applied on it. The running time and corresponding relative error of algorithms are averaged over 50 Monte Carlo simulations. To show the acceleration of randomized algorithms against the deterministic counterparts, we first report the running time of the THOSVD and the STHOSVD algorithms and their randomized variants. The results for the noiseless data tensor are displayed in Figure \ref{FIGrand1} and the results for the noisy case (20 dB) are displayed in Figure \ref{FIGrand2}. For computing the leading singular vectors of unfolding matrices, we used both the truncated SVD (MATLAB function svds) and the EVD (MATLAB function eigs) and compared their executive times. In displayed figures, the notations ``svds'' and ``eigs'' mean that the underlying algorithms exploit the truncated SVD and the EVD for computing the factor matrices, respectively. Our results show that the EVD is much faster than using the svds for computing the factor matrices and their randomized versions make the algorithms much faster. The performance comparison of deterministic and randomized algorithms for the noiseless case are reported in Table \ref{Table1}. It is seen that randomized algorithms provide approximately the same results as the deterministic ones but with less computational time. For the noisy data tensor with 20 dB, the fit of the THOSVD, the STHOSVD algorithms and their randomized variants were almost the same and equal to $90\%$ but the randomized algorithms were much faster as seen in Figure \ref{FIGrand2}. For the noiseless data tensor, we did not use any power iteration or oversampling parameter. For the noisy case, we set power iteration $q=2$ and oversampling $P=10$, and this is why for a data tensor with the same size and multilinear rank, more acceleration is achieved using randomized algorithms for the noiseless data tensor compared to the noisy data tensor.

\begin{figure}[!h]
\includegraphics[width=1\linewidth]{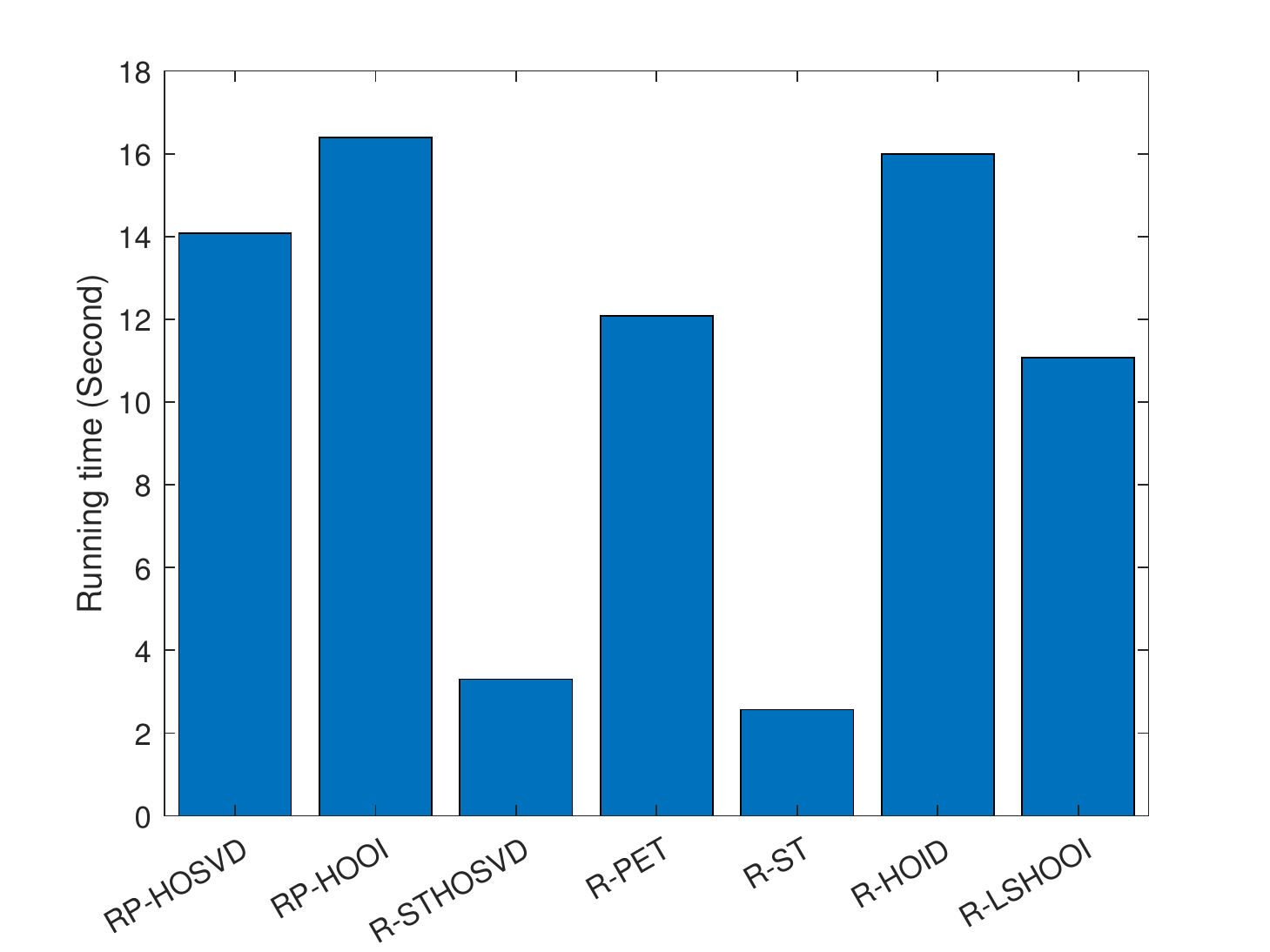}
\centering
\caption{\small Running time comparison of various randomized algorithms for a noiseless random data tensor of size $1000\times 1000\times 1000$ and multilinear rank $(20,40,30)$.}
\label{FIGComRan}
\end{figure}

We also compared the performance of different randomized algorithms for the noiseless data tensor in Figure \ref{FIGComRan} (also roughly the same results were obtained for the noisy data tensor and we skip them). The results show that the R-STHOSVD algorithm (Algorithm \ref{RSeqHOSVD}) and the randomized sampling Tucker (R-ST) algorithm (Algorithm \ref{SamTucker}) are the fastest algorithms for computation of the HOSVD or the Tucker decomposition. Note the R-ST algorithm based on uniform sampling is also relatively fast but when the noise level is relatively high or the singular values of unfolding matrices do not decay very fast, its accuracy is degraded due to additive error bounds. In such situations, the leverage-score probability distribution can provide more accurate solutions with relative error bound guarantee but requiring higher computational cost. 

%
%

\subsubsection{\bf Function Based Tensors}
\begin{figure}[!h]
\includegraphics[width=1\linewidth]{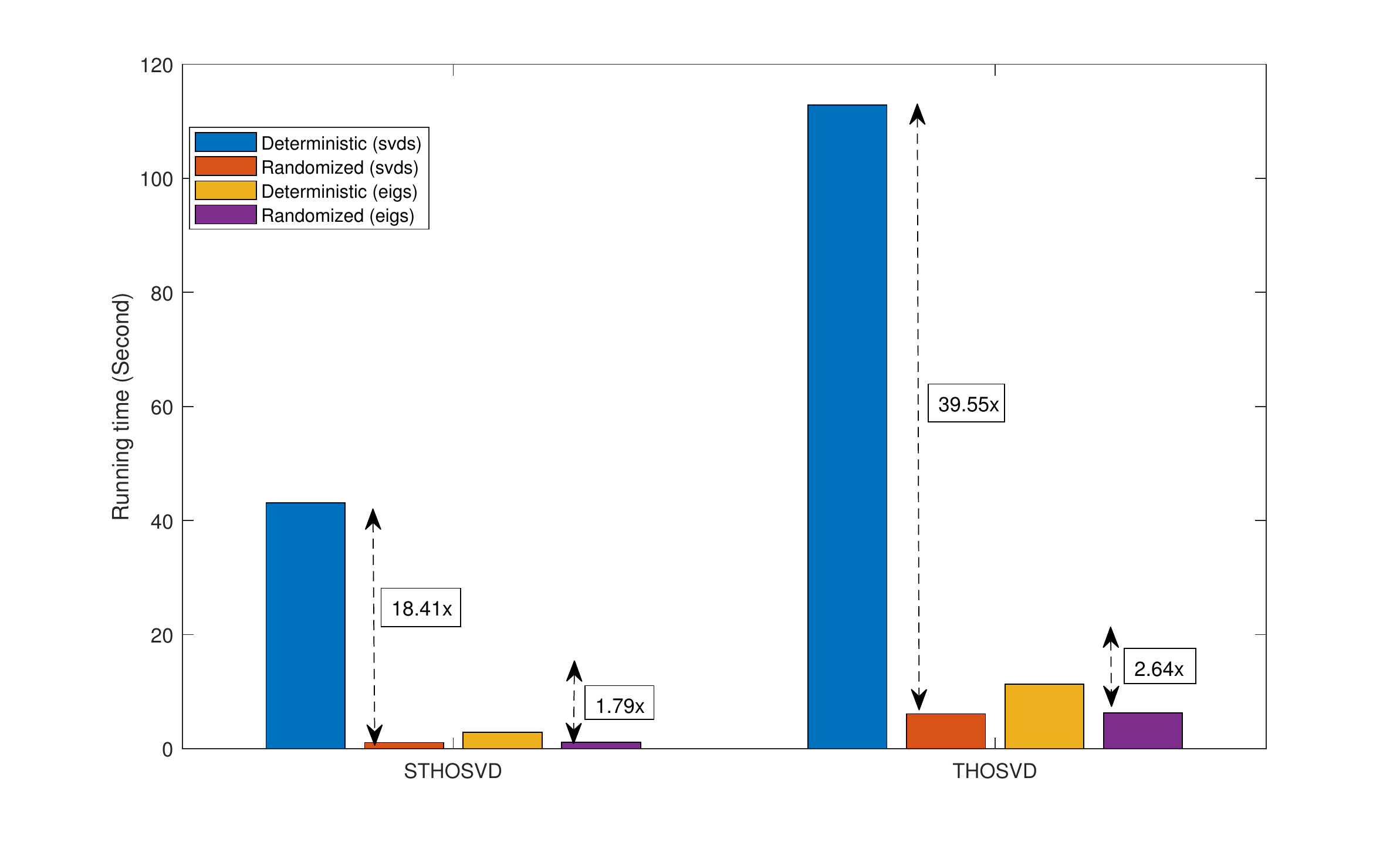}
\centering
\caption{\small Running time comparison of deterministic and randomized variants of the THOSVD and the STHOSVD algorithms (svds and eigs versions) for the function based tensor \eqref{funtensor} of size $1000\times 1000\times 1000$ and multilinear rank $(30,30,30)$.}
\label{FIGFunBased}
\end{figure}

\begin{figure}[!h]
\includegraphics[width=1\linewidth]{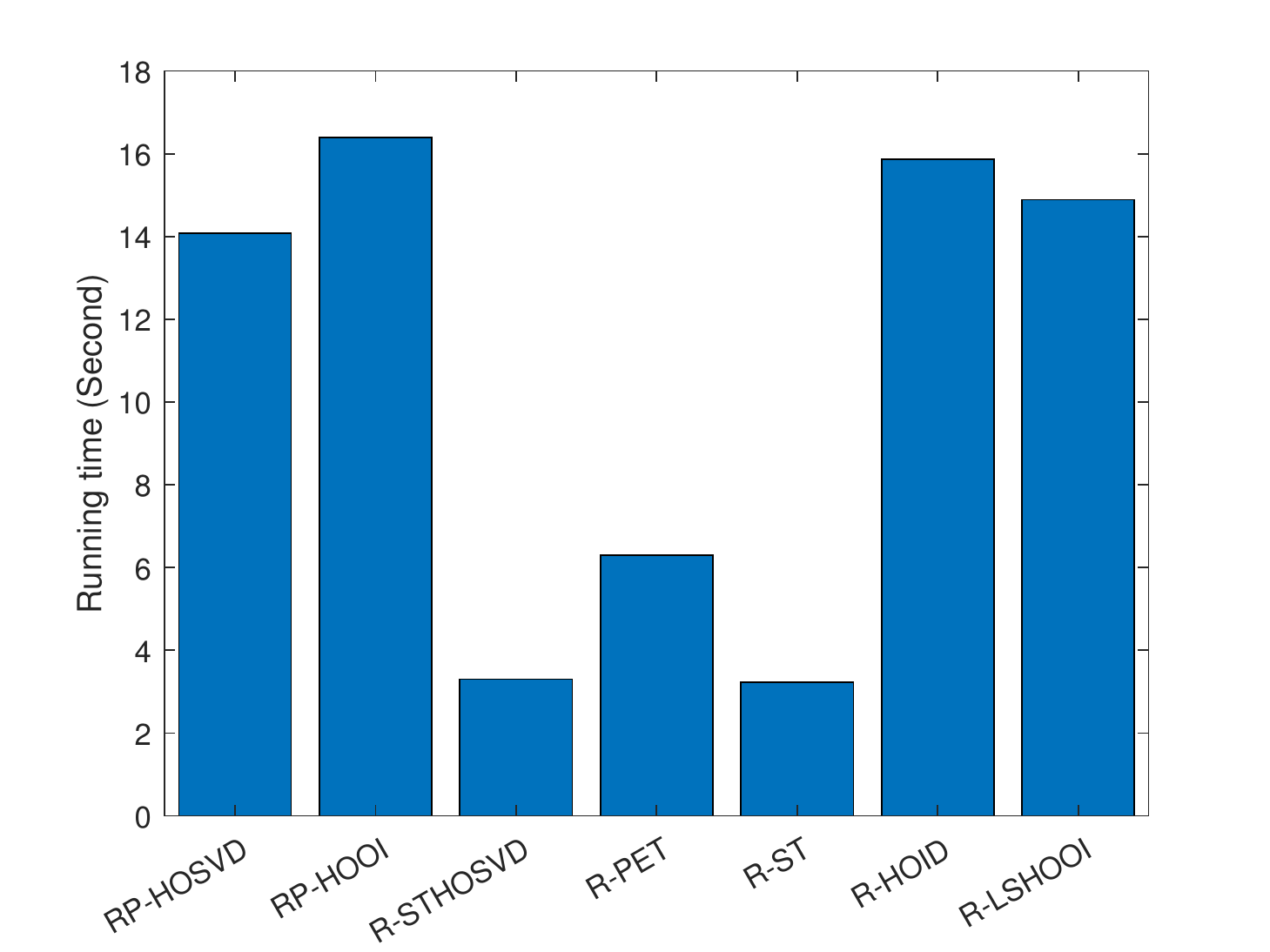}
\centering
\caption{\small Running time comparison of various randomized algorithms for the function based tensor \eqref{funtensor} of size $1000\times 1000\times 1000$ and multilinear rank $(30,30,30)$.}
\label{FIGFunBased2}
\end{figure}

Consider a 3rd order tensor of size $1000\times 1000\times 1000$ whose components are generated as follows
\begin{equation}\label{funtensor}
\underline{\mathbf X}(i,j,k)=\frac{1}{{\sqrt[5]{{{i^5} + {j^5} + {k^5}}}}}.
\end{equation}
It is shown in \cite{beylkin2002numerical, beylkin2005algorithms, chinnamsetty2007tensor, grasedyck2004existence, khoromskij2008tensor, hackbusch2005hierarchical, tyrtyshnikov2003tensor} that such a data tensor admits a very low multilinear rank approximation since the singular values of the corresponding unfolding matrices decay very fast, (for more related funcation based tensors, see \cite{beylkin2002numerical, beylkin2005algorithms, chinnamsetty2007tensor, grasedyck2004existence, khoromskij2008tensor, hackbusch2005hierarchical, tyrtyshnikov2003tensor} and the references therein). We set the multilinear rank as $(30,30,30)$ in our computations and apply deterministic and randomized algorithms on the mentioned data tensor. The running time of the THOSVD, the STHOSVD and their randomized variants are reported in Figure \ref{FIGFunBased}. Also, the performance comparison of deterministic and randomized algorithms are reported in Table \ref{Table1}. Here again the simulation results indicate the superiority of randomized algorithms over the deterministic ones. To compare the running time of different randomized algorithms, see Figure \ref{FIGFunBased2}. The results show that the R-STHOSVD and the R-ST algorithms are the most efficient randomized algorithms for the low multilinear rank approximation of the data tensor \eqref{funtensor}.

\tabcolsep=0.11cm
 \begin{table*}
 \centering
   \caption{Relative error comparison of different deterministic and randomized algorithms for computing low mutilinear rank approximation of data tensors \eqref{noiselesstensor}, \eqref{funtensor}, \eqref{SPT}  and \eqref{HilbertTensor}.}
   \footnotesize
 \begin{tabular}{||c c c c c c c c c c c||}
 \hline
 Data & HOOI & STHOSVD & HOSVD & RP-HOSVD & RP-HOOI & R-STHOSVD & R-PET & R-ST & R-HOID & R-LSHOOI\\ [0.5ex] 
  \hline\hline
Data tensor \eqref{noiselesstensor} & 1.0303e-13 & 1.5103e-13 & 1.5203e-13 & 6.0334e-13 &  3.0012e-13 &  6.0334e-13 & 1.5620e-12 & 1.2812e-11 & 4.7531e-13 & 6.2167e-13 \\ 
 \hline\hline
Data tensor \eqref{funtensor} & 1.0031e-8 & 1.6095e-8 & 1.6124e-8 & 3.6206e-8 &  1.0031e-8 &  3.1637e-8& 3.4651e-7 & 4.1674e-6 & 6.2213e-8 & 3.4352e-8 \\ 
 \hline
 Sparse tensor \eqref{SPT}  & 2.7732e-6 & 4.8484e-6 & 9.5498e-6 & 9.5498e-6&2.7732e-6 & 7.5750e-6 & 7.5750e-5 & 7.5750e-4 & 7.5750e-6 & 4.7531e-6 \\
 \hline
 Hilbert tensor \eqref{HilbertTensor}  & 1.0044e-10 & 4.9682e-10 & 8.4250e-10 & 1.1041e-10 &  1.2122e-10 & 1.1478e-10 & 4.2234e-9 & 6.3753e-8 & 6.9878e-10 & 4.7531e-10\\
 \hline
 \end{tabular}\label{Table1}
\end{table*}

\subsubsection{\bf Random Sparse Tensors}
\begin{figure}[!h]
\includegraphics[width=1\linewidth]{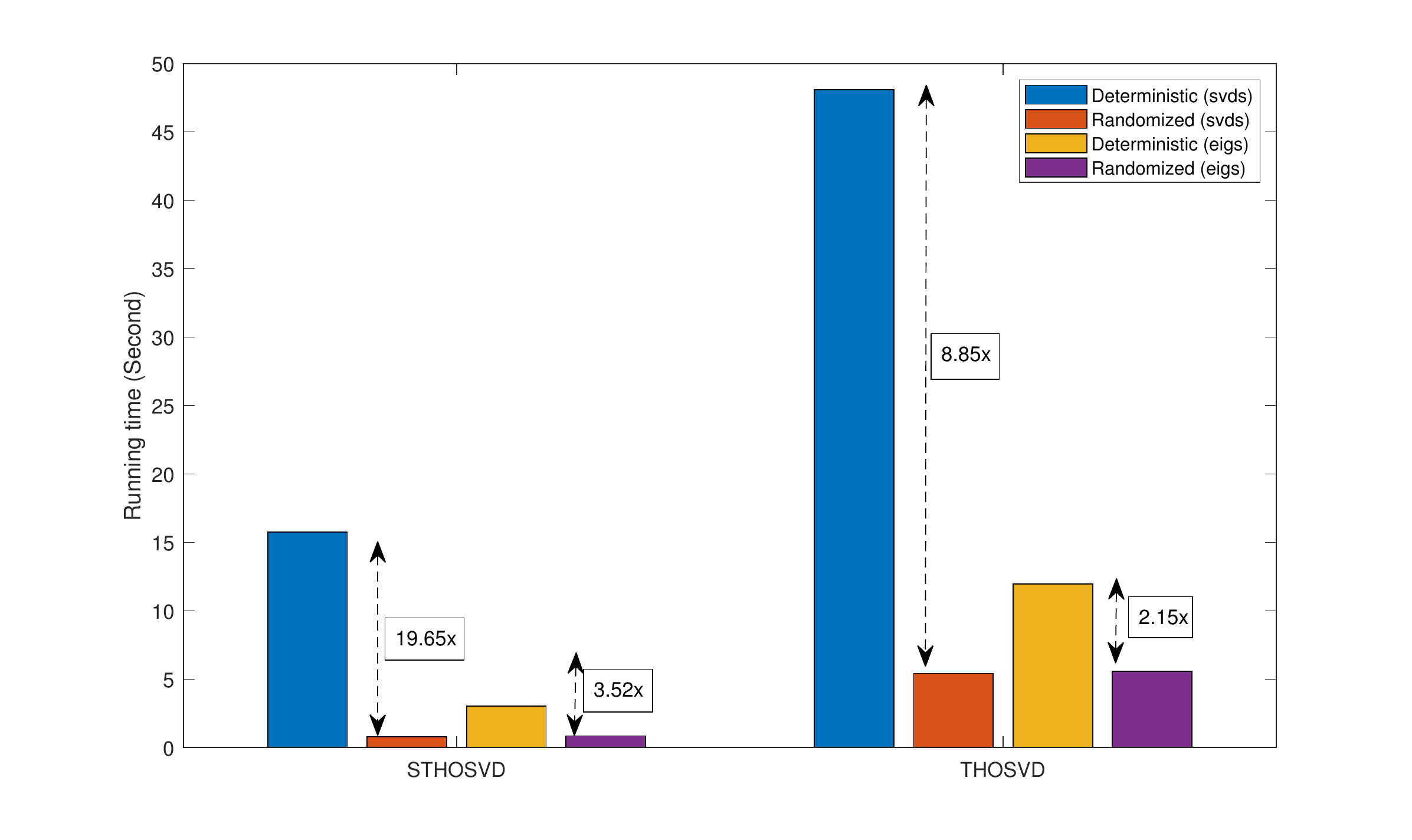}
\centering
\caption{\small Running time comparison of deterministic and randomized variants of the THOSVD and the STHOSVD algorithms (svds and eigs versions) for the sparse tensor \eqref{SPT} of size $1000\times 1000\times 1000$ and multilinear rank $(15,15,15)$.}
\label{sparse}
\end{figure}

\begin{figure}[!h]
\includegraphics[width=1\linewidth]{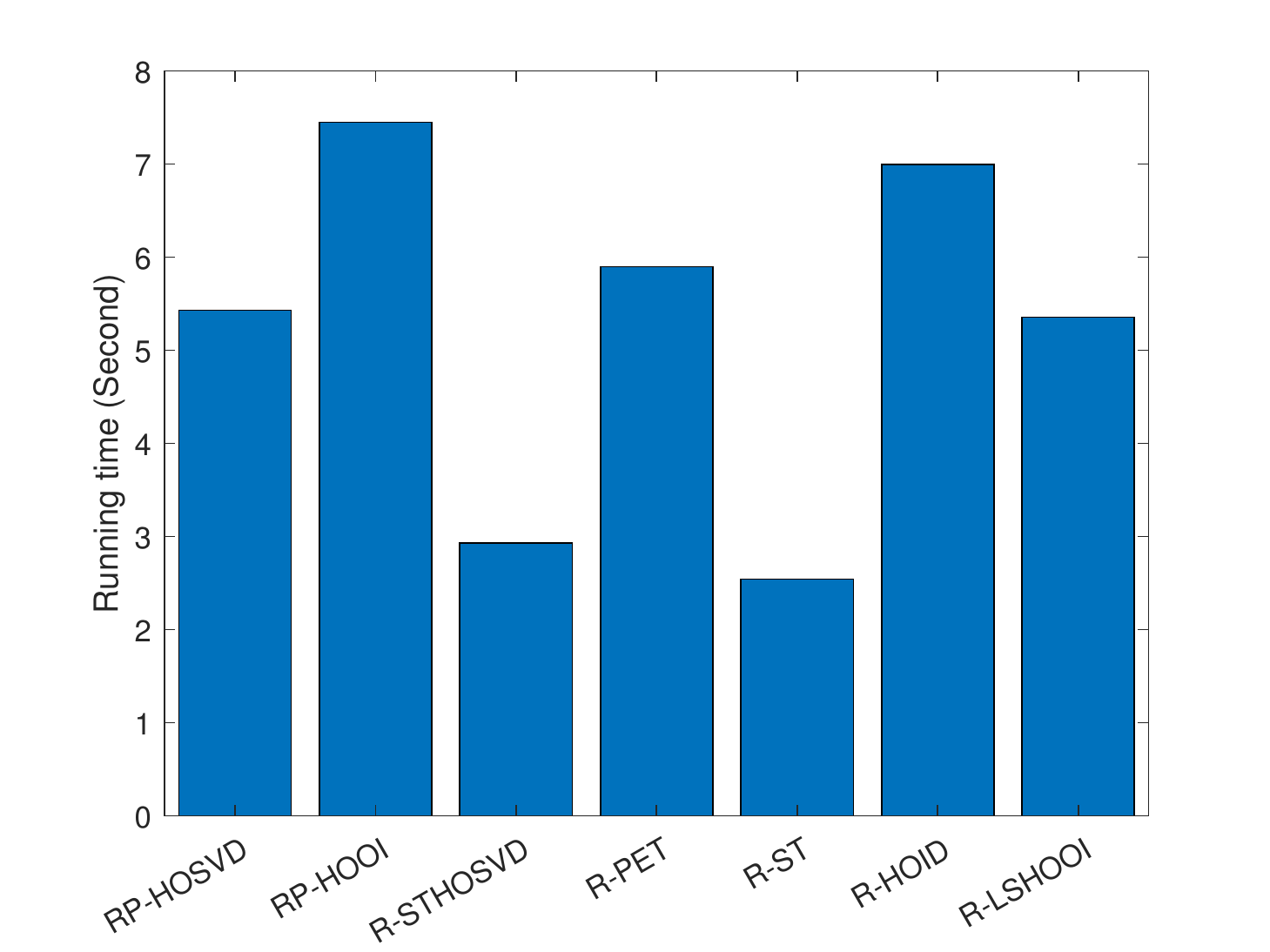}
\centering
\caption{\small Running time comparison of various randomized algorithms for the sparse tensor \eqref{SPT} of size $1000\times 1000\times 1000$ and multilinear rank $(15,15,15)$.}
\label{sparse2}
\end{figure}

Consider the following data tensor
\begin{equation}\label{SPT}
\underline{\mathbf X} = \sum\limits_{i = 1}^{10} {\frac{\gamma }{{{i^2}}}{{\mathbf x}_i} \circ {{\mathbf y}_i} \circ {{\mathbf z}_i}}  + \sum\limits_{i = 11}^{200} {\frac{1 }{{{i^2}}}{{\mathbf x}_i} \circ {{\mathbf y}_i} \circ {{\mathbf z}_i}}, 
\end{equation}
where ${\mathbf x}_i,\,{\mathbf y}_i$ and ${\mathbf z}_i$ are sparse vectors with only $0.05$ sparsity (percent  non-zero numbers). We set $I_1=I_2=I_3=1000$ and $\gamma=1000$ and generate a 3rd-order tensor $\underline{\mathbf X}$.
The parameters used in this simulation are the same as previous experiments and we use a multilinear rank $R=(15,15,15)$. The running time of the THOSVD and the STHOSVD algorithms and their randomized forms are reported in Figure \ref{sparse}. Also, the performance comparison of deterministic and randomized algorithms are reported in Table \ref{Table1}. Figure \ref{sparse} and Table \ref{Table1}, show that the randomized algorithms can significantly speed-up the decomposition procedure with achieving roughly the same accuracy as the deterministic algorithms. In Figure \ref{sparse2}, we make a comparison among different randomized algorithms. The results show that the R-STHOSVD and the R-ST algorithms are the most promising approaches for the low multilinear approximation of the data tensor \eqref{SPT}. It is also seen that the randomized algorithms can achieve almost the same accuracy as the deterministic ones but they are much faster.

\subsubsection{\bf Hilbert Tensors}
\begin{figure}[!h]
\includegraphics[width=1\linewidth]{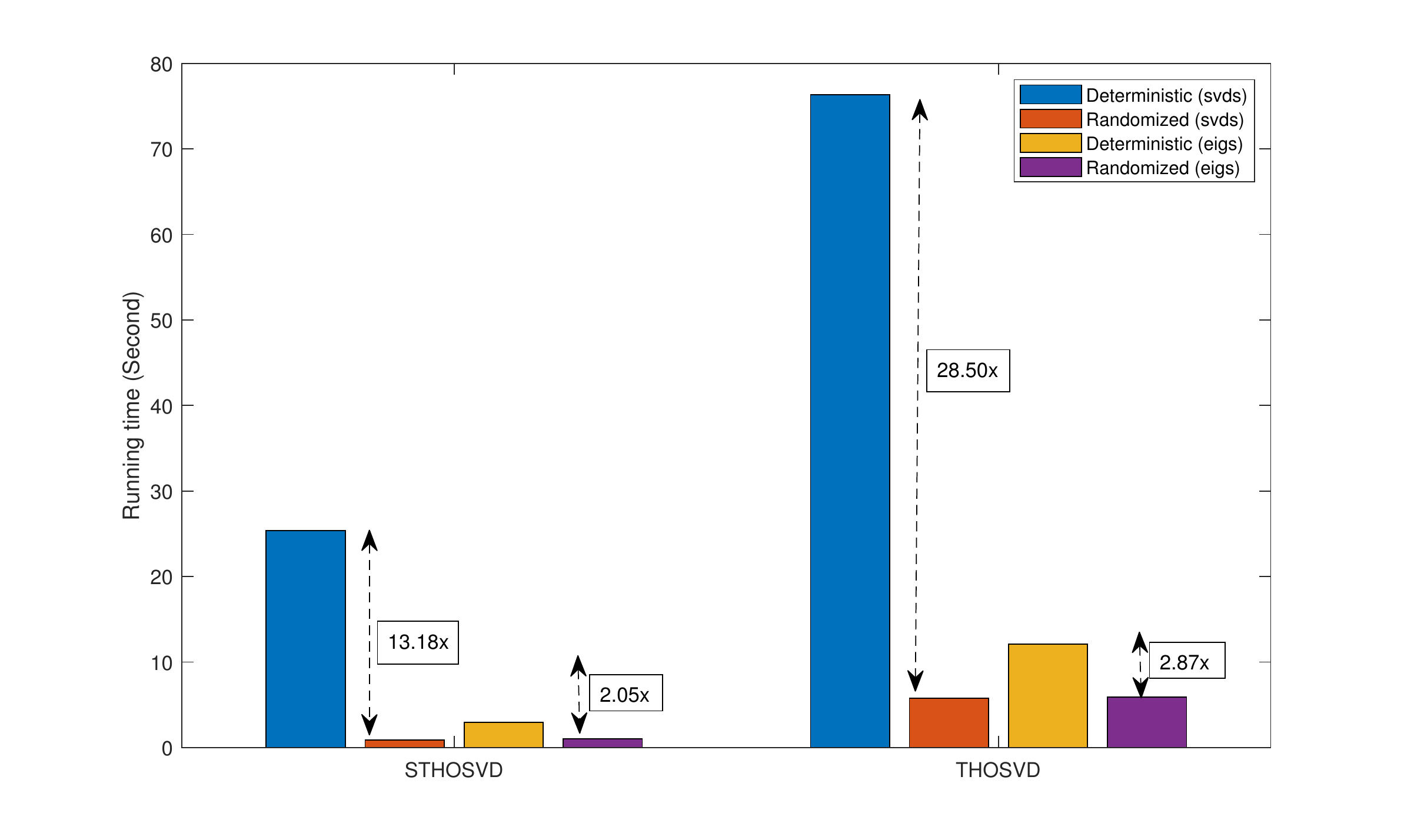}
\centering
\caption{\small Running time comparison of deterministic and randomized variants of the THOSVD and the STHOSVD algorithms (svds and eigs versions) for the Hilbert tensor of size $1000\times 1000\times 1000$ and multilinear rank $(20,20,20)$.}
\label{Hilbert}
\end{figure}

\begin{figure}[!h]
\includegraphics[width=1\linewidth]{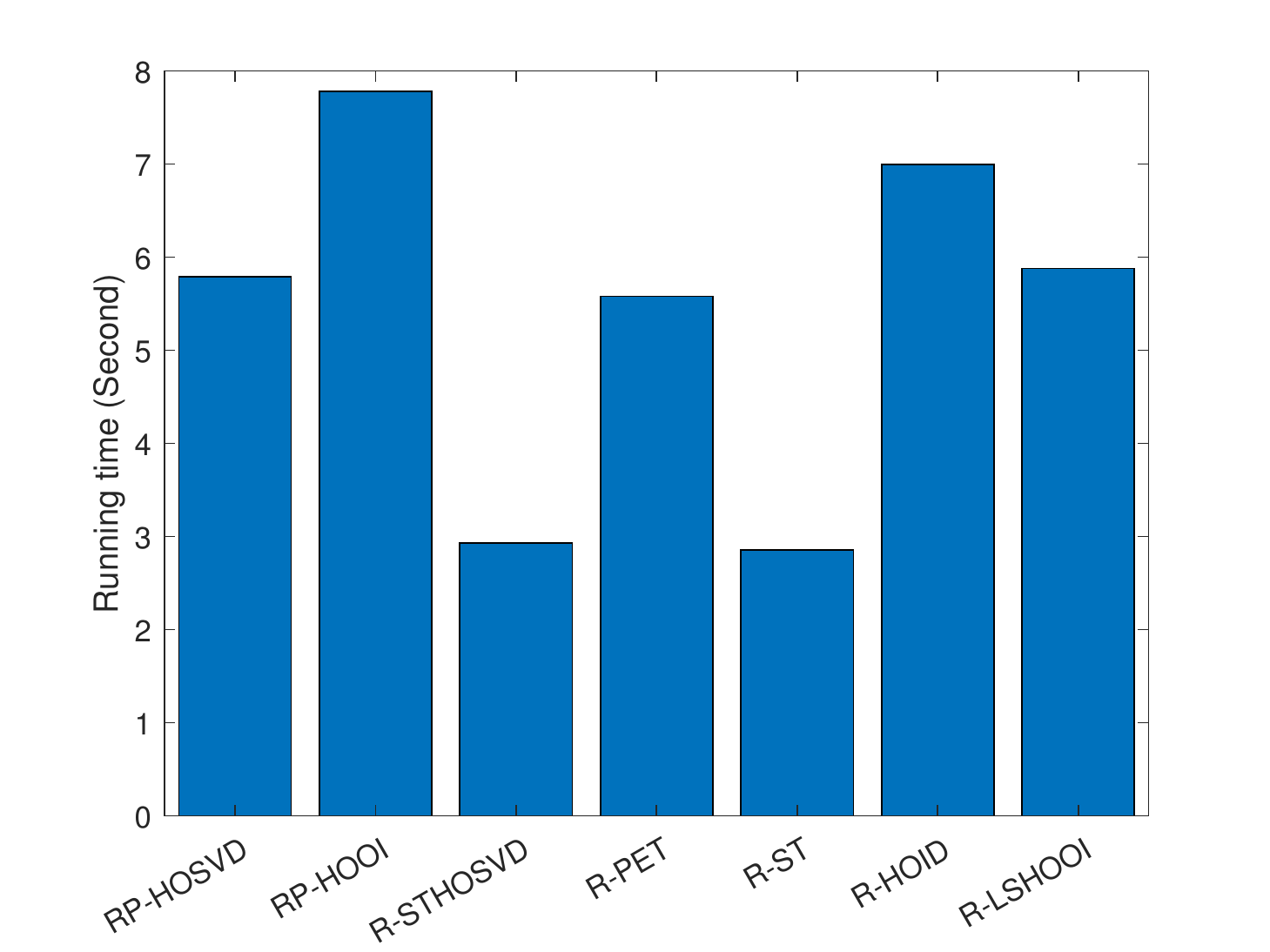}
\centering
\caption{\small Running time comparison of various randomized algorithms for the Hilbert tensor of size $1000\times 1000\times 1000$ and multilinear rank $(20,20,20)$.}
\label{Hilbert2}
\end{figure}

An $N$th order Hilbert tensor is defined as follows
\begin{equation}\label{HilbertTensor}
\underline{\mathbf X}(i_1,i_2,\ldots,i_N)=\frac{1}{i_1+i_2+\ldots+i_N-N+1},
\end{equation}
which is a natural generalization of Hilbert matrices (second order tensors). We generate a 3th-order Hilbert tensor of size $1000\times 1000\times 1000$ and use approximate multilinear rank $(20,20,20)$ within the randomized algorithms. The running time of the THOSVD and the STHOSVD algorithms and their randomized variants are reported in Figure \ref{Hilbert}. Also, the performance comparison of deterministic and randomized algorithms are reported in Table \ref{Table1}. Once again the superiority of randomized algorithms over the deterministic ones is visible. In Figure \ref{Hilbert2}, we make a comparison among different randomized algorithms. The results show that the R-STHOSVD and the R-ST algorithms are the best and most efficient algorithms for the low multilinear rank approximation of the data tensor \eqref{HilbertTensor}. It is seen that the randomized algorithms can achieve almost the same accuracy as the deterministic ones but much faster.

\subsection{Test 2 -- real data tensors} 

\subsubsection{COIL-100 data set}
\begin{figure}[!h]
\includegraphics[width=1\linewidth]{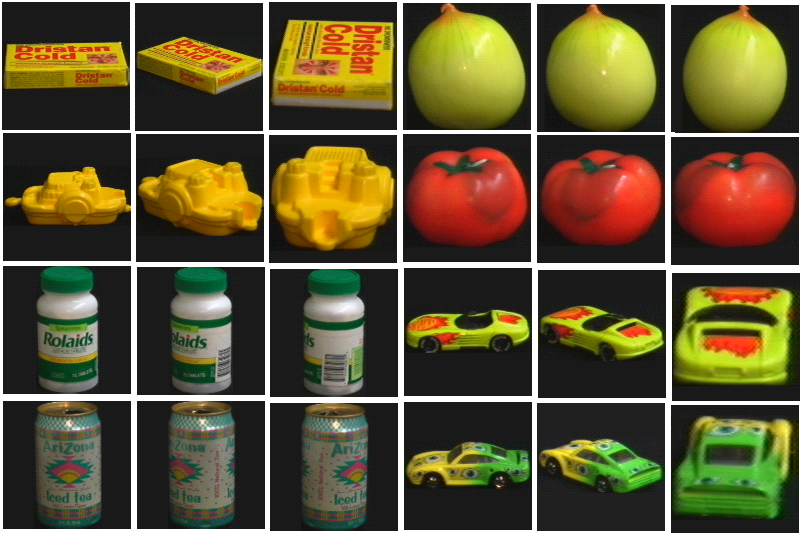}
\centering
\caption{\small The first 8 images of the COIL-100 dataset for three random rotations.}
\label{coilsamples}
\end{figure}

\begin{figure}[!h]
\includegraphics[width=1\linewidth]{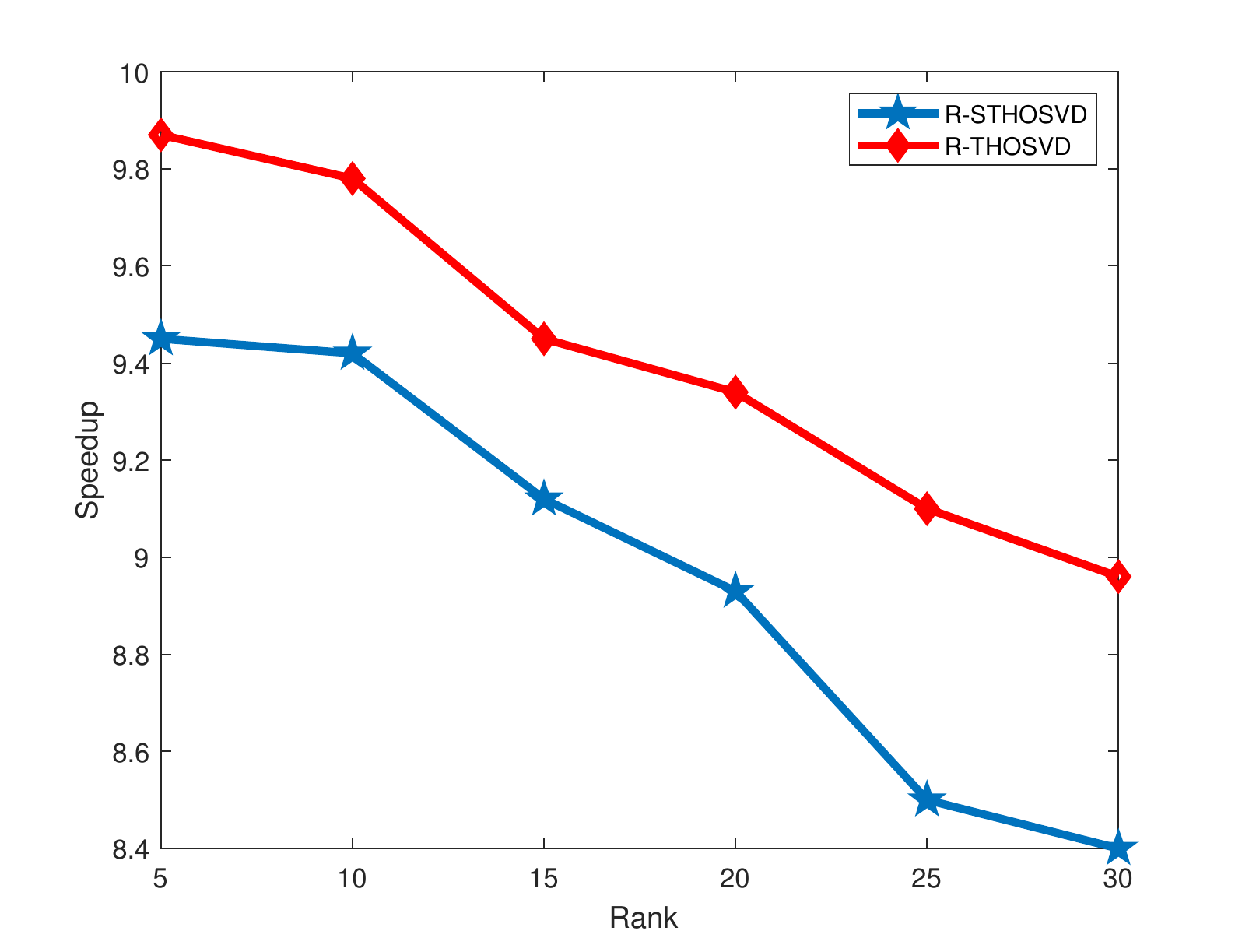}
\centering
\caption{\small The running time speed-up of randomized STHODV and randomized THOSVD for decomposing COIL-100 for selected multilinear rank $(R,R,R)$.}
\label{speedupSHOSVD}
\end{figure}

\begin{figure}[!h]
\includegraphics[width=1\linewidth]{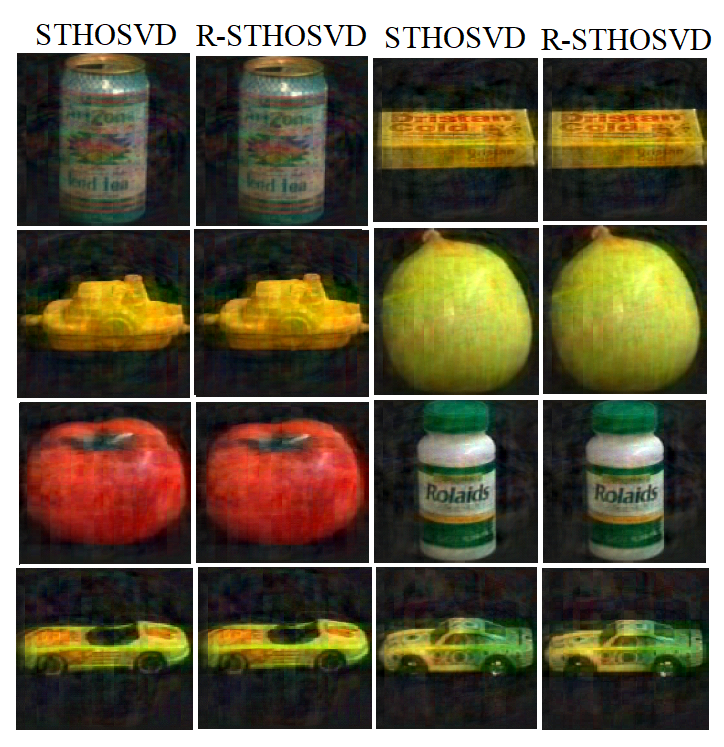}
\centering
\caption{\small The reconstructed images using randomized and deterministic STHOSVD algorithms for the first 8 images.}
\label{Reconstruction}
\end{figure}

Columbia object image library COIL-100 is a dataset which consists of 7200 color images (100 objects under 72 different rotations) \cite{nene1996columbia}. Each image is of size $128\times 128\times 3$ and the whole dataset naturally can be represented as a 5th-order tensor of size $128\times 128\times 3\times 72\times 100$. The first 5 images and some random rotations are depicted in Figure \ref{coilsamples}. We reshape this tensor to a 3rd-order tensor of size $1024\times 1152\times 300$ and apply THSOVD and STHOSVD algorithms and their randomized  variants on this 3rd order data tensor. The speedup achieved by the R-STHOSVD and the R-THOSVD for different uniform multilinear ranks $(R,R,R),\,\,R=5,10,15,\ldots,30$ are shown in Figure \ref{speedupSHOSVD}. Also, the reconstructed images using the STHOSVD and R-STHOSVD for multilinear rank $(450,450,250)$ are displayed in Figure \ref{Reconstruction}. The results indicate that the randomized algorithms provide roughly the same performance as the deterministic algorithms but in much less running times.

\begin{figure}[!h]
\includegraphics[width=1\linewidth]{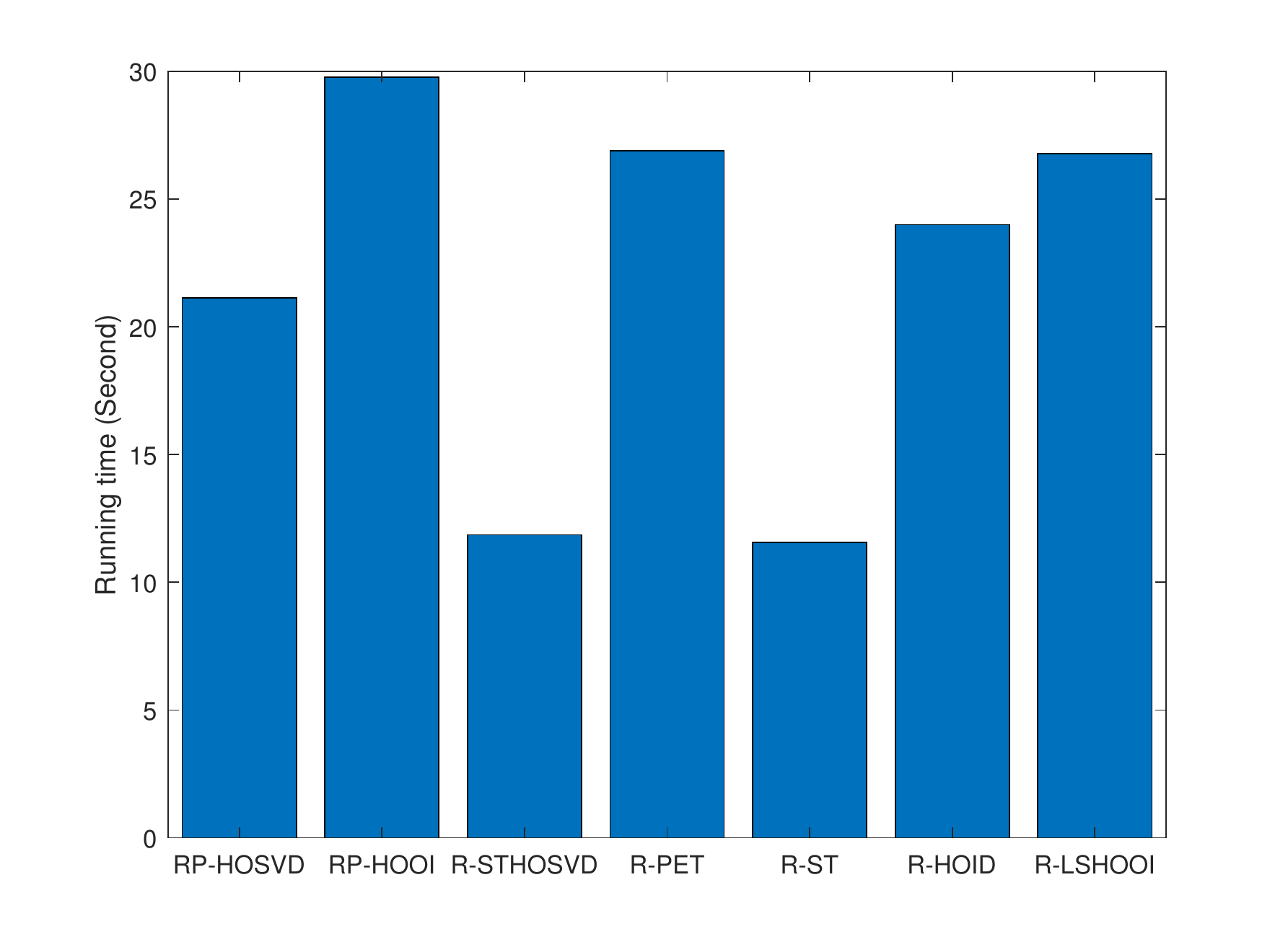}
\centering
\caption{\small The running time comparison of various randomized algorithms for decomposing COIL100 data set and multilinear rank $(450,450,250)$.}
\label{coil100compare}
\end{figure}

The running time comparison of various randomized algorithms for multilinear rank $(20,20,20)$ are reported in Figure \ref{coil100compare}. Our results again show the effectiveness of the R-STHOSVD for decomposing large-scale data tensors.

%
%
%

\subsubsection{NELL-2 data tensor (Large-scale sparse tensor)}
In this experiment, we show the applicability and effectiveness of the randomized algorithms for handling large-scale sparse data tensors. To this end, we considered the {\it NELL-2} \cite{carlson2010toward} which is accessible at the formidable repository of open sparse tensors and tool (FROSTT) database \cite{frosttdataset}. It is a 3rd order tensor of size $12092 \times 9184 \times 28818$ with 76879419 number of nonzero elements. The NELL-2 dataset was tested on machine learning algorithms which exploit relationship between variables\cite{carlson2010toward}. It was generated at Chicago university as a part of ``Read the Web" project.

For this sparse data tensor, all the deterministic algorithms introduced in the paper were not able to decompose the mentioned data tensor in the Tucker format due to high memory and computational complexity issues while the randomized STHOSVD was able to decompose it. The running time and relative-error achieved by the RSTHOSVD algorithm for the uniorm multilinear rank $(R,R,R),\,R=100,125,150\ldots,200$ are reported in Table \ref{TableRelativeRunning}.

\subsection{video datasets} 
\begin{figure*}
\centering
     \includegraphics[width=.5\linewidth]{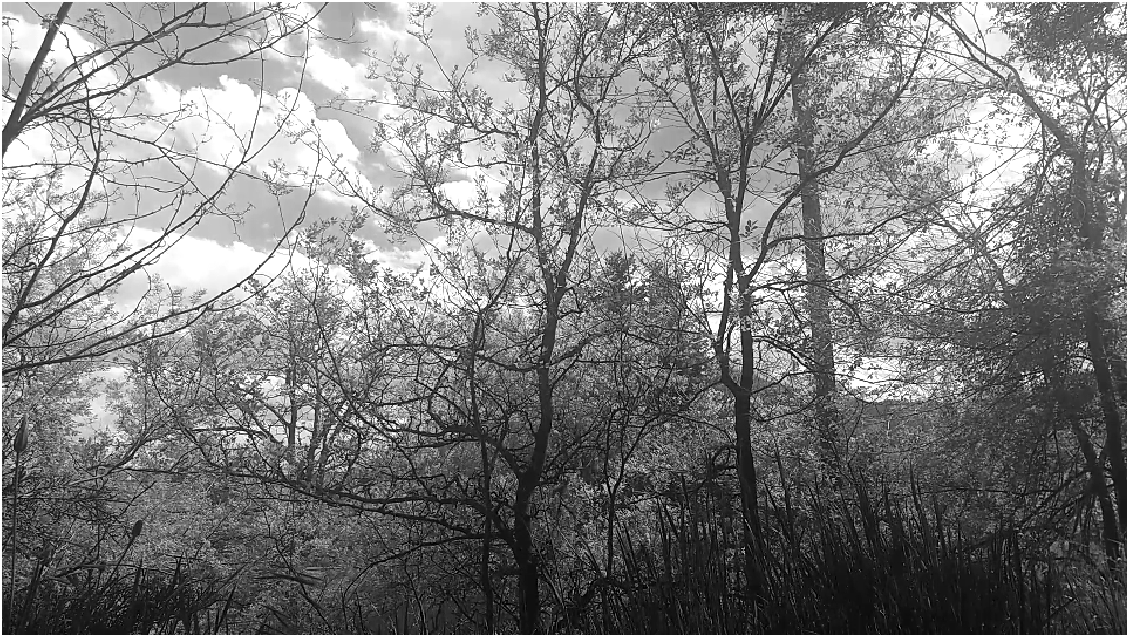}
    \caption{Visualization of the original frame number 50 of the grayscale video.}
    \label{fig:videooriginalframe}
\end{figure*}

\begin{figure*}
{
          \includegraphics[width=1\linewidth]{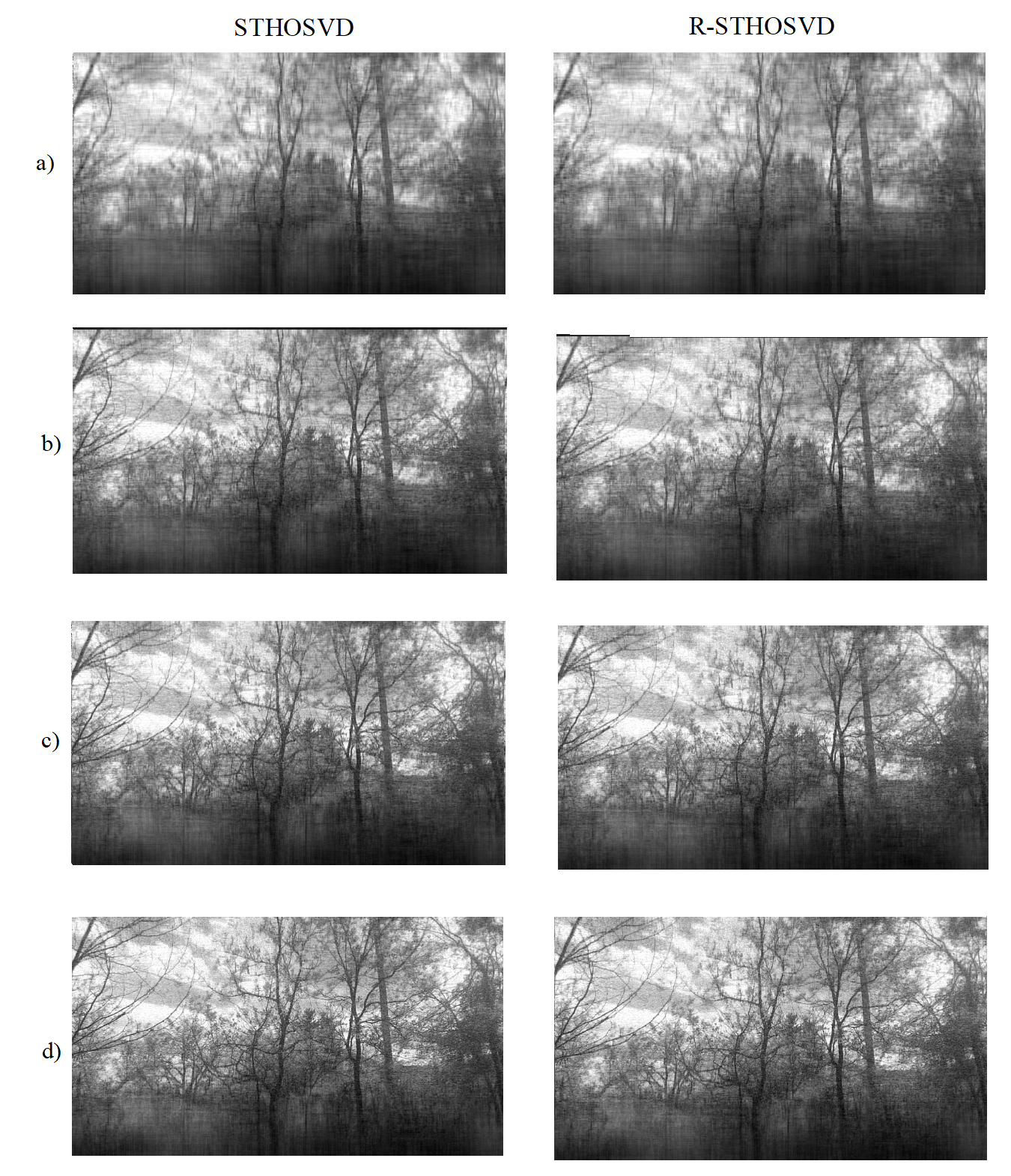}
    }
    \caption{Visualization of the reconstructed frame number 50 using deterministic and randomized STHOSVD algorithms for different multilinear ranks, a) Multilinear rank $(50,50,50),\,\,\,$ b) Multilinear rank $(100,100,100),\,\,\,$ c) Multilinear rank $(150,150,150,)\,\,\,$ d) Multilinear rank $(200,200,200)$.}
    \label{fig:framereconstruction}
\end{figure*}

In this experiment we compress a large-scale grayscale video \href{https://drive.google.com/file/d/1usBNBSfnPuy1S2Oy8-QQPusrvBmVTohl/view} in the Tucker format. The data tensor was firstly used to classify the frames in \cite{malik2018low}. The original data tensor is a 4rd order tensor of size $2200 \times 1080\times 3 \times 1980$ consisting of $1980$ RGB frames. For our task, we consider its grayscale form which is of size $2200 \times 1080\times 1980$. The frame number 50 is displayed in Figure \ref{fig:videooriginalframe}. 
We applied the THOSVD and the STHOSVD algorithms on the mentioned grayscale video for different uniform multilinear rank $(R,R,R),\,R=10,20,\ldots,200$. The relative error and running time comparison are reported in Table \ref{TableRelativeRunningError}. Also, the reconstructed images using randomized and deterministic STHOSVD algorithms are visualized in Figure \ref{fig:framereconstruction}. In both scenarios, it is seen that the randomized algorithms achieve roughly the same relative error as the deterministic algorithms but much faster.




\begin{table*}[h!]
\centering
\caption{The running time (in second) and relative errors of the randomized STHOSVD Algorithm for NELL-2 dataset.}\label{TableRelativeRunning}
\begin{tabular}{c|c|c|c|c}
\hline
&\multicolumn{2}{c}{STHOSVD}&\multicolumn{2}{c}{R-STHOSVD}\\
\hline
Target Multilinear Rank & Running time & Relative error& Running time & Relative error\\
\hline
(50,50,50) &  \xmark &   \xmark &   212.56 &   0.6248\\
(75,75,75) &   \xmark &  \xmark & 442.72  & 0.5720   \\
(100,100,100) &   \xmark &  \xmark &   853.25&    0.5255\\
(125,125,125) &   \xmark &  \xmark &  1347.94 & 0.4944 \\
(150,150,150) &   \xmark &  \xmark &   1412.56 &   0.4695\\
(175,175,175) &   \xmark &  \xmark &  1514.43 & 0.4484  \\
(200,200,200) &   \xmark &  \xmark &   2009.51 &   0.4283\\
\hline
\end{tabular}
\end{table*}

\begin{table*}[h!]
\centering
\caption{The relative errors of deterministic and randomized algorithms for the grayscale video.}\label{TableRelativeRunningError}
\begin{tabular}{c|c|c|c|c|||c|c|c|c|c}
\hline
&\multicolumn{4}{c}{Relatiev error}&\multicolumn{4}{c}{Running time (Second)}\\
\hline
Target multilinear rank &\multicolumn{1}{c}{HOSVD}&\multicolumn{1}{c}{R-HOSVD}&\multicolumn{1}{c}{STHOSVD}&\multicolumn{1}{c}{R-STHOSVD}&\multicolumn{1}{c}{HOSVD}&\multicolumn{1}{c}{R-HOSVD}&\multicolumn{1}{c}{STHOSVD}&\multicolumn{1}{c}{R-STHOSVD}&1/Compression ratio\\
\hline
(10,10,10) & 0.3226 & 0.3346 & 0.3210 & 0.3247  & 122.90 & 17.42 & 38.42 & 3.40 &  1.1393e-05 \\
(20,20,20) & 0.3030 & 0.3178 & 0.3020 & 0.3069  & 140.11 & 19.98 & 41.20 & 3.82 & 2.4062e-05 \\
(30,30,30) & 0.2913 & 0.3069 & 0.2901 & 0.2975  & 133.59 & 21.29 & 36.09 & 4.17  & 3.9282e-05 \\
(40,40,40) & 0.2827 & 0.2989 & 0.2814 & 0.2894  & 136.55 & 24.60 & 42.04 & 5.59  & 5.8327e-05 \\
(50,50,50) & 0.2761 & 0.2937 & 0.2742 & 0.2833  & 133.60 & 29.81 & 43.49 & 6.85  & 8.2475e-05 \\
(60,60,60) & 0.2699 & 0.2879 & 0.2681 & 0.2772  & 158.87 & 34.81 & 48.91 & 7.90  & 1.1300e-04 \\
(70,70,70) & 0.2647 & 0.2835 & 0.2625 & 0.2728  & 159.74 & 38.46 & 49.35 & 8.77  & 1.5118e-04 \\
(80,80,80) & 0.2600 & 0.2779 & 0.2575 & 0.2678  & 157.28 & 39.18 & 50.22 & 10.51  & 1.9828e-04 \\
(90,90,90) & 0.2553 & 0.2740 & 0.2528 & 0.2634  & 156.29 & 44.98 & 50.34 & 10.26  & 2.5559e-04 \\
(100,100,100) & 0.2509 & 0.2697 & 0.2485 & 0.2597  & 158.16 & 47.69 & 52.07 & 13.15&3.2437e-04\\
(110,110,110) & 0.2471 & 0.2656 & 0.2445 & 0.2559 & 159.13 & 53.78 & 50.61 & 12.94 &4.0591e-04\\
(120,120,120) & 0.2432 & 0.2618 & 0.2407 & 0.2522  & 156.14 & 54.72 & 51.77 & 15.01 & 5.0148e-04\\
(130,130,130) & 0.2397 & 0.2583 & 0.2371 & 0.2486 & 157.94 & 58.61 & 53.65 & 16.41  & 6.1235e-04\\
(140,140,140) & 0.2364 & 0.2544 & 0.2337 & 0.2451  & 159.34 & 56.99 & 53.01 & 17.66 & 7.3981e-04\\
(150,150,150) & 0.2333 & 0.2511 & 0.2305 & 0.2421 & 163.63 & 65.13 & 55.07 & 17.65 & 8.8511e-04\\
(160,160,160) & 0.2301 & 0.2475 & 0.2273 & 0.2389 & 160.37 & 68.05 & 54.04 & 20.17 & 0.0010 \\
(170,170,170) & 0.2270 & 0.2443 & 0.2243 & 0.2359  & 161.34 & 70.54 & 54.93 & 21.15 & 0.0012\\
(180,180,180) & 0.2239 & 0.2408 & 0.2213 & 0.2331  & 159.83 & 74.41 & 56.63 & 22.73 & 0.0014\\
(190,190,190) & 0.2209 & 0.2378 & 0.2184 & 0.2301 & 161.92 & 72.89 & 57.45 & 23.21  & 0.0017\\
(200,200,200) & 0.2181 & 0.2346 & 0.2156 & 0.2271 & 162.20 & 82.63 & 57.64 & 25.59  & 0.0019\\
\hline
\end{tabular}
\end{table*}

%

\section{Conclusion}\label{Sec:Concl}
In this paper, we reviewed and extended a variety of state-of-the-art randomized algorithms for computing the Tucker decomposition and the Higher Order SVD (HOSVD). We studied both single-pass and multi-pass randomized algorithms and also random projection and sampling techniques for computing the Tucker decomposition and the HOSVD. Simulations were conducted on synthetic and real datasets to compare the performance and running time of the deterministic and the randomized algorithms with particular focus on showing the superiority of randomized algorithms over the deterministic ones.

\section*{Acknowledgments}
The first author wishes to thank Dr. Varvara Logacheva for her constructive comments and useful discussion on the paper. This work was partially supported by the Ministry of Education and Science of the Russian Federation (grant project (14.756.31.0001).

\ifCLASSOPTIONcaptionsoff
  \newpage
\fi
\bibliographystyle{ieeetr}

      \bibliography{test}
 \EOD

\end{document}